\documentclass[aos,MSNbibl,seceqn,citesort,dvips]{arximspdf}
\usepackage{graphicx}

\doi{10.1214/11-AOS949}
\volume{40}
\issue{4}
\pubyear{2012}
\firstpage{1935}
\lastpage{1967}

\makeatletter

\newcommand{\OH}{O\,H}

\newcommand{\tr}{\operatorname{tr}}
\newcommand{\lh}{\ell}
\newcommand{\F}{\mathcal{I}}
\newcommand{\s}{\mathcal{S}}
\newcommand{\lr}{\mathcal{L}}
\newcommand{\degr}{\operatorname{deg}}
\newcommand{\ts}{\mathcal{Y}}
\newcommand{\add}{\mathcal{A}}
\newcommand{\proj}{\mathcal{P}}
\newcommand{\g}{g_\gamma}
\newcommand{\C}{\mathcal{C}}
\newcommand{\Sa}{S^\ast}
\newcommand{\La}{L^\ast}

\newtheorem{theorem}{Theorem}[section]
\newtheorem{corl}[theorem]{Corollary}
\newtheorem{lemm}[theorem]{Lemma}
\newtheorem{prop}[theorem]{Proposition}

\newproclaim{defin}{Definition}[section]

\makeatother

\begin{document}
\begin{frontmatter}

\title{Latent variable graphical model selection via convex optimization\thanksref{T1}}
\runtitle{Latent variable model selection}
\relateddoi{}{Discussed in \doi{10.1214/12-AOS979},
\doi{10.1214/12-AOS980}, \doi{10.1214/12-AOS981},
\doi{10.1214/12-AOS984}, \doi{10.1214/12-AOS985} and
\doi{10.1214/12-AOS1001};
rejoinder at \doi{10.1214/12-AOS1020}.}

\thankstext{T1}{Supported in part by AFOSR Grant FA9550-08-1-0180,
in part under a MURI through AFOSR Grant FA9550-06-1-0324,
in part under a MURI through AFOSR Grant FA9550-06-1-0303
and in part by NSF FRG 0757207.}

\begin{aug}
\author[A]{\fnms{Venkat}~\snm{Chandrasekaran}\corref{}\ead[label=e1]{venkatc@caltech.edu}},
\author[B]{\fnms{Pablo~A.}~\snm{Parrilo}\ead[label=e2]{parrilo@mit.edu}}
\and\\
\author[B]{\fnms{Alan~S.}~\snm{Willsky}\ead[label=e3]{willsky@mit.edu}}
\runauthor{V. Chandrasekaran, P. A. Parrilo and A. S. Willsky}
\affiliation{California Institute of Technology,
Massachusetts Institute of Technology and~Massachusetts Institute of Technology}
\address[A]{V. Chandrasekaran\\
Department of Computing\\
\quad and Mathematical Sciences\\
California Institute of Technology\\
Pasadena, California 91106\\
USA\\
\printead{e1}}
\address[B]{P. A. Parrilo\\
A. S. Willsky\\
Laboratory for Information\\
\quad and Decision Systems\\
Department of Electrical Engineering\\
\quad and Computer Science\\
Massachusetts Institute of Technology\\
Cambridge, Massachusetts 02139\\
USA\\
\printead{e2}\\
\phantom{E-mail: }\printead*{e3}} 
\end{aug}

\received{\smonth{8} \syear{2010}}
\revised{\smonth{11} \syear{2011}}

%
\begin{abstract}
Suppose we observe samples of a \textit{subset} of a collection of random
variables. No additional information is provided about the number of
latent variables, nor of the relationship between the latent and
observed variables. Is it possible to discover the number of latent
components, and to learn a statistical model over the entire collection
of variables? We address this question in the setting in which the
latent and observed variables are jointly Gaussian, with the
conditional statistics of the observed variables conditioned on the
latent variables being specified by a graphical model. As a first step
we give natural conditions under which such latent-variable Gaussian
graphical models are identifiable given marginal statistics of only the
observed variables. Essentially these conditions require that the
conditional graphical model among the observed variables is sparse,
while the effect of the latent variables is ``spread out'' over most of
the observed variables. Next we propose a tractable convex program
based on regularized maximum-likelihood for model selection in this
latent-variable setting; the regularizer uses both the $\ell_1$ norm
and the nuclear norm. Our modeling framework can be viewed as a
combination of dimensionality reduction (to identify latent variables)
and graphical modeling (to capture remaining statistical structure not
attributable to the latent variables), and it consistently estimates
both the number of latent components and the conditional graphical
model structure among the observed variables. These results are
applicable in the high-dimensional setting in which the number of
latent/observed variables grows with the number of samples of the
observed variables. The geometric properties of the algebraic
varieties of sparse matrices and of low-rank matrices play an important
role in our analysis.
\end{abstract}

%

\end{frontmatter}

\section{Introduction and setup}
\label{secintro}

Statistical model selection in the high-dimen\-sional regime arises in a
number of applications. In many data analysis problems in geophysics,
radiology, genetics, climate studies, and image processing, the number
of samples available is comparable to or even smaller than the number
of variables. As empirical statistics in these settings may not be
well-behaved (see \mbox{\cite{MarP1967,Joh2001}}), high-dimensional model
selection is therefore both challenging and of great interest. A model
selection problem that has received considerable attention recently is
the estimation of covariance matrices in the high-dimensional setting.
As the sample covariance matrix is poorly behaved in such a regime,
some form of \textit{regularization} of the sample covariance is adopted
based on assumptions about the true underlying covariance matrix \cite
{BicL2008a,Elk2008,BicL2008b,LedW2003,WuP2003,FanFL2008}.


\subsection*{Graphical models} A number of papers have studied
covariance estimation in the context of \textit{Gaussian graphical model
selection}. A \textit{Gaussian graphical model}~\cite{SpeK1986,Lau1996}
(also commonly referred to as a Gauss--Markov random field) is a
statistical model defined with respect to a graph, in which the nodes
index a collection of jointly Gaussian random variables and the edges
represent the conditional independence relations (Markov structure)
among the variables. In such models the sparsity pattern of the inverse
of the covariance matrix, or the \textit{concentration} matrix, directly
corresponds to the graphical model structure. Specifically, consider a
Gaussian graphical model in which the covariance matrix is given by a
positive-definite $\Sigma^\ast$ and the concentration matrix is given
by $K^\ast= (\Sigma^\ast)^{-1}$. Then an edge $\{i,j\}$ is present in
the underlying graphical model if and only if $K^\ast_{i,j} \neq0$. In
particular the absence of an edge between two nodes implies that the
corresponding variables are independent conditioned on all the other
variables. The model selection method usually studied in such a
Gaussian graphical model setting is $\ell_1$-regularized
maximum-likelihood, with the $\ell_1$ penalty applied to the entries of
the concentration matrix to induce sparsity. The consistency properties
of such an estimator have been studied \cite
{RotBLZ2008,RavWRY2008,LamF2009}, and under suitable conditions \cite
{LamF2009,RavWRY2008} this estimator is also ``sparsistent,'' that is,
the estimated concentration matrix has the same sparsity pattern as the
true model from which the samples are generated. An alternative
approach to $\ell_1$-regularized maximum-likelihood is to estimate the
sparsity pattern of the concentration matrix by performing regression
separately on each variable~\cite{MeiB2006}; while such a method
consistently estimates the sparsity pattern, it does not directly
provide estimates of the covariance or concentration matrix.

In many applications throughout science and engineering (e.g.,
psychology, computational biology, and economics), a challenge is that
one may not have access to observations of all the relevant phenomena,
that is, some of the relevant variables may be latent or unobserved. In
general latent variables pose a significant difficulty for model
selection because one may not know the number of relevant latent
variables, nor the relationship between these variables and the
observed variables. Typical algorithmic methods that try to get around
this difficulty usually fix the number of latent variables as well as
the structural relationship between latent and observed variables
(e.g., the graphical model structure between latent and observed
variables), and use the EM algorithm to fit parameters \cite
{DemLR1977}. This approach suffers from the problem that one optimizes
nonconvex functions, and thus one may get stuck in suboptimal local
minima. An alternative suggestion~\cite{EliNN2007} is one based on a
greedy, local, combinatorial heuristic that assigns latent variables to
groups of observed variables, via some form of clustering of the
observed variables; however, this approach has no consistency guarantees.

\subsection*{Our setup} In this paper we study the problem of
latent-variable graphical model selection in the setting where all the
variables, both observed and latent, are jointly Gaussian. More
concretely, $X$ is a Gaussian random vector in $\mathbb{R}^{p+h}$, $O$
and $H$ are disjoint subsets of indices in $\{1,\ldots,p+h\}$ of
cardinalities $|O| = p$ and $|H| = h$, and the corresponding subvectors
of $X$ are denoted by $X_O$ and $X_H$, respectively. Let the covariance
matrix underlying $X$ be denoted by $\Sigma^\ast_{(\OH)}$. The
marginal statistics corresponding to the observed variables $X_O$ are
given by the marginal covariance matrix $\Sigma^\ast_O$, which is
simply a submatrix of the full covariance matrix $\Sigma^\ast_{(\OH
)}$. However, suppose that we parameterize our model by the
concentration matrix $K^\ast_{(\OH)} = (\Sigma^\ast_{(\OH)})^{-1}$,
which as discussed above reveals the connection to graphical models.
Here the submatrices $K^\ast_O, K^\ast_{O,H}, K^\ast_H$ specify (in the
full model) the dependencies among the observed variables, between the
observed and latent variables, and among the latent variables,
respectively. In such a parameterization, the \textit{marginal
concentration matrix} $(\Sigma^\ast_O)^{-1}$ corresponding to the
observed variables $X_O$ is given by the Schur complement \cite
{HorJ1990} with respect to the block $K^\ast_H$:
%
\begin{equation} \label{eqschur}
\tilde{K}^\ast_{O} = (\Sigma^\ast_O)^{-1} = K^\ast_O - K^\ast_{O,H}
(K^\ast_H)^{-1} K^\ast_{H,O}.
\end{equation}
Thus if we only observe the variables $X_O$,
we only have access to $\Sigma^\ast_O$ (or~$\tilde{K}^\ast_O$). The two
terms that compose $\tilde{K}^\ast_O$ above have interesting
properties. The matrix $K^\ast_O$ specifies the concentration matrix of
the \textit{conditional statistics} of the observed variables given the
latent variables. If these conditional statistics are given by a sparse
graphical model, then $K^\ast_O$ is \textit{sparse}. On the other hand,
the matrix $K^\ast_{O,H} (K^\ast_H)^{-1} K^\ast_{H,O}$ serves as a
\textit
{summary} of the effect of marginalization over the latent variables
$X_H$. This matrix has small rank if the number of latent, unobserved
variables $X_H$ is small relative to the number of observed variables
$X_O$. Therefore the marginal concentration matrix $\tilde{K}^\ast_O$
is generally \textit{not sparse} due to the additional low-rank term
$K^\ast_{O,H} (K^\ast_H)^{-1} K^\ast_{H,O}$. Hence standard graphical
model selection techniques applied directly to the observed variables
$X_O$ are not useful.

A modeling paradigm that infers the effect of the latent variables
$X_H$ would be more suitable in order to provide a concise explanation
of the underlying statistical structure. Hence we approximate the
sample covariance by a model in which the concentration matrix \textit
{decomposes} into the sum of a sparse matrix and a low-rank matrix,
which reveals the conditional graphical model structure in the observed
variables as well as the number of and effect due to the unobserved
latent variables. Such a method can be viewed as a blend of principal
component analysis and graphical modeling. In standard graphical
modeling one would directly approximate a concentration matrix by a
sparse matrix to learn a sparse graphical model, while in principal
component analysis the goal is to explain the statistical structure
underlying a set of observations using a small number of latent
variables (i.e., approximate a covariance matrix as a low-rank matrix).
In our framework we learn a sparse graphical model among the observed
variables \textit{conditioned} on a few (additional) latent variables.
These latent variables are \textit{not} principal components, as the
conditional statistics (conditioned on these latent variables) are
given by a graphical model. Therefore we refer to these latent
variables informally as \textit{latent components}.

\subsection*{Contributions} Our first contribution in Section \ref
{seciden} is to address the fundamental question of \textit
{identifiability} of such latent-variable graphical models given the
marginal statistics of only the observed variables. The critical point
is that we need to tease apart the correlations induced due to
marginalization over the latent variables from the conditional
graphical model structure among the observed variables. As the
identifiability problem is one of \textit{uniquely} decomposing the sum
of a sparse matrix and a low-rank matrix into the individual
components, we study the algebraic varieties of sparse matrices and
low-rank matrices. An important theme in this paper is the connection
between the tangent spaces to these algebraic varieties and the
question of identifiability. Specifically let $\Omega(K^\ast_O)$ denote
the tangent space at $K^\ast_O$ to the algebraic variety of sparse
matrices, and let $T(K^\ast_{O,H} (K^\ast_H)^{-1} K^\ast_{H,O})$ denote
the tangent space at $K^\ast_{O,H} (K^\ast_H)^{-1} K^\ast_{H,O}$ to the
algebraic variety of low-rank matrices. Then the \textit{statistical}
question of identifiability of $K^\ast_O$ and $K^\ast_{O,H} (K^\ast
_H)^{-1} K^\ast_{H,O}$ given $\tilde{K}^\ast_O$ is determined by the
\textit{geometric} notion of \textit{transversality} of the tangent spaces
$\Omega(K^\ast_O)$ and $T(K^\ast_{O,H} (K^\ast_H)^{-1} K^\ast_{H,O})$.
The study of the transversality of these tangent spaces leads to
natural conditions for identifiability. In particular we show that
latent-variable models in which $(1)$ the sparse matrix $K^\ast_O$ has
a small number of nonzeros per row/column, and $(2)$ the low-rank
matrix $K^\ast_{O,H} (K^\ast_H)^{-1} K^\ast_{H,O}$ has row/column
spaces that are not closely aligned with the coordinate axes, are
identifiable. These conditions have natural statistical
interpretations. The first condition ensures that there are no densely
connected subgraphs in the conditional graphical model structure among
the observed variables, that is, that these conditional statistics are
indeed specified by a sparse graphical model. Such statistical
relationships may otherwise be mistakenly attributed to the effect of
marginalization over some latent variable. The second condition ensures
that the effect of marginalization over the latent variables is
``spread out'' over many observed variables; thus, the effect of
marginalization over a latent variable is not confused with the
conditional graphical model structure among the observed variables. In
fact the first condition is often assumed in standard graphical model
selection without latent variables (e.g.,~\cite{RavWRY2008}).


As our next contribution we propose a \textit{regularized
maximum-likelihood decomposition} framework to approximate a given
sample covariance matrix by a model in which the concentration matrix
decomposes into a sparse matrix and a low-rank matrix. Based on the
effectiveness of the $\ell_1$ norm as a tractable convex relaxation for
recovering sparse models~\cite{Don2006a,Don2006b,CanRT2006} and the
nuclear norm for low-rank matrices~\cite{Faz2002,RecFP2009,CanR2009},
we propose the following penalized likelihood method given a sample
covariance matrix $\Sigma^n_O$ formed from $n$ samples of the observed
variables:
%
\begin{eqnarray}\label{eqsdp}
(\hat{S}_n,\hat{L}_n) &=& \mathop{\arg\min}_{S,L}  -\lh(S-L; \Sigma
^n_O)
+ \lambda_n  \bigl(\gamma\|S\|_{1} +
\tr(L)\bigr)
\nonumber\\[-8pt]\\[-8pt]
&&\hspace*{19pt}\mbox{s.t. } S-L \succ0,  L \succeq0.\nonumber
\end{eqnarray}
The constraints $\succ0$ and $\succeq0$ impose positive-definiteness
and positive-semi\-definiteness. The function $\lh$ represents the
Gaussian log-likelihood $\lh(K;\Sigma) = \log\det(K) - \tr(K
\Sigma)$
for $K \succ0$, where $\tr$ is the trace of a matrix and $\det$ is the
determinant. Here $\hat{S}_n$ provides an estimate of $K^\ast_O$, which
represents the conditional concentration matrix of the observed
variables; $\hat{L}_n$ provides an estimate of $K^\ast_{O,H} (K^\ast
_H)^{-1} K^\ast_{H,O}$, which represents the effect of marginalization
over the latent variables. The regularizer is a combination of the
$\ell
_1$ norm applied to $S$ and the nuclear norm applied to $L$ (the
nuclear norm reduces to the trace over the cone of symmetric,
positive-semidefinite matrices), with $\gamma$ providing a trade-off
between the two terms. This variational formulation is a \textit{convex
optimization} problem, and it is a regularized max-det program that can
be solved in polynomial time using general-purpose solvers~\cite{WanST2009}.

Our main result in Section~\ref{secmain} is a proof of the consistency
of the estimator~(\ref{eqsdp}) in the high-dimensional regime in which
both the number of observed variables and the number of latent
components are allowed to grow with the number of samples (of the
observed variables). We show that for a suitable choice of the
regularization parameter $\lambda_n$, there exists a range of values of
$\gamma$ for which the estimates $(\hat{S}_n,\hat{L}_n)$ have the same
sparsity (and sign) pattern and rank as $(K^\ast_O,K^\ast_{O,H}
(K^\ast
_H)^{-1} K^\ast_{H,O})$ with high probability (see Theorem \ref
{theomain}). The key technical requirement is an identifiability
condition for the two components of the marginal concentration matrix
$\tilde{K}^\ast_O$ with respect to the Fisher information (see
Section~\ref{subsecfi}). We make connections between our condition and
the irrepresentability conditions required for support/graphical-model
recovery using $\ell_1$ regularization \cite
{ZhaY2006,RavWRY2008,Wai2009}. Our results provide numerous scaling
regimes under which consistency holds in latent-variable graphical
model selection. For example, we show that \textit{under suitable
identifiability conditions consistent model selection is possible even
when the number of samples and the number of latent variables are on
the same order as the number of observed variables} (\textit{see Section \ref
{subsecscal}}).

\subsection*{Related previous work} The problem of decomposing the sum
of a sparse matrix and a low-rank matrix via convex optimization into
the individual components was initially studied in~\cite{ChaSPW2009} by
a superset of the authors of the present paper, with conditions derived
under which the convex program exactly recovers the underlying
components. In subsequent work Cand\`es et al.~\cite{CanLMW2009} also
studied this sparse-plus-low-rank decomposition problem, and provided
guarantees for exact recovery using the convex program proposed in
\cite
{ChaSPW2009}. The problem setup considered in the present paper is
quite different and is more challenging because we are only given
access to an inexact sample covariance matrix, and we wish to produce
an \textit{inverse} covariance matrix that can be decomposed as the sum
of sparse and low-rank components (preserving the sparsity pattern and
rank of the components in the true underlying model). In addition to
proving the consistency of the estimator~(\ref{eqsdp}), we also
provide a statistical interpretation of our identifiability conditions
and describe natural classes of latent-variable Gaussian graphical
models that satisfy these conditions. As such our paper is closer in
spirit to the many recent papers on covariance selection, but with the
important difference that some of the variables are not observed.

\subsection*{Outline} Section~\ref{secbg} gives some background and a
formal problem statement. Section~\ref{seciden} discusses the
identifiability question, Section~\ref{secmain} states the main
results of this paper, and Section~\ref{secproofs} gives some proofs.
We provide experimental demonstration of the effectiveness of our
estimator on synthetic and real data in Section~\ref{secsims}, and
conclude with a brief discussion in Section~\ref{secconc}. Some of our
technical results are deferred to supplementary material~\cite{supp}.


\section{Problem statement and background}
\label{secbg}

We give a formal statement of the latent-variable model selection
problem. We also briefly describe various properties of the algebraic
varieties of sparse matrices and of low-rank matrices, and the
properties of the Gaussian likelihood function.

The following matrix norms are employed throughout this paper. $\|M\|
_2$ denotes the spectral norm, or the largest singular value of $M$. $\|
M\|_\infty$ denotes the largest entry in magnitude of $M$. $\|M\|_F$
denotes the Frobenius norm, or the square root of the sum of the
squares of the entries of $M$. $\|M\|_\ast$ denotes the nuclear norm,
or the sum of the singular values of $M$ (this reduces to the trace for
positive-semidefinite matrices). $\|M\|_1$ denotes the sum of the
absolute values of the entries of $M$. A~number of \textit{matrix
operator} norms are also used. For example, let $\mathcal{Z}\dvtx \mathbb
{R}^{p \times p} \rightarrow\mathbb{R}^{p \times p}$ be a linear
operator acting on matrices. Then the induced operator norm is defined
as $\|\mathcal{Z}\|_{q \rightarrow q} \triangleq{\max_{N \in\mathbb
{R}^{p \times p}, \|N\|_q \leq1}} \|\mathcal{Z}(N)\|_q$.
Therefore,\vspace*{1pt}
$\|\mathcal{Z}\|_{F \rightarrow F}$ denotes the spectral norm of the
operator $\mathcal{Z}$. The only vector norm used is the Euclidean
norm, which is denoted by $\| \cdot\|$. Given any norm \mbox{$\|\cdot\|_q$}
(either a vector norm, a matrix norm or a matrix operator norm), the
dual norm is given by $\|M\|_q^\ast\triangleq\operatorname{sup}\{\langle M,
N \rangle | \|N\|_q \leq1\}$.

\subsection{Problem statement}
\label{subsecps}

In order to analyze latent-variable model selection methods, we need to
define an appropriate notion of model selection consistency for
latent-variable graphical models. Given the two components $K^\ast_O$
and $K^\ast_{O,H} (K^\ast_H)^{-1} K^\ast_{H,O}$ of the concentration
matrix of the marginal distribution~(\ref{eqschur}), there are \textit
{infinitely} many configurations of the latent variables [i.e.,
matrices $K^\ast_H \succ0, K^\ast_{O,H} = (K^\ast_{H,O})^T$] that give
rise to the \textit{same} low-rank matrix $K^\ast_{O,H} (K^\ast_H)^{-1}
K^\ast_{H,O}$. Specifically for any nonsingular matrix $B \in\mathbb
{R}^{|H| \times|H|}$, one can apply the transformations $K^\ast_H
\rightarrow B K^\ast_H B^T, K^\ast_{O,H} \rightarrow K^\ast_{O,H} B^T$
and still preserve the low-rank matrix $K^\ast_{O,H} (K^\ast_H)^{-1}
K^\ast_{H,O}$. In \textit{all} of these models the marginal statistics of
the observed variables $X_O$ remain the same upon marginalization over
the latent variables $X_H$. The key \textit{invariant} is the low-rank
matrix $K^\ast_{O,H} (K^\ast_H)^{-1} K^\ast_{H,O}$, which summarizes
the effect of marginalization over the latent variables. Consequently,
from here on we use the notation $\Sa= K^\ast_O$ and $\La= K^\ast
_{O,H} (K^\ast_H)^{-1} K^\ast_{H,O}$. These observations give rise to
the following notion of structure recovery.
%
\begin{defin}
A pair of $|O| \times|O|$ symmetric matrices $(\hat{S},\hat{L})$ is an
\textit{algebraically correct} estimate of a latent-variable Gaussian
graphical model given by the concentration matrix $K^\ast_{(\OH)}$ if
the following conditions hold:
\begin{longlist}[(3)]
\item[(1)] The sign-pattern of $\hat{S}$ is the same as that of $\Sa$ [here
$\operatorname{sign}(0) = 0$]:
\[
\operatorname{sign}(\hat{S}_{i,j}) =
\operatorname{sign}(\Sa_{i,j})\qquad
\forall i,j.
\]

\item[(2)] The rank of $\hat{L}$ is the same as the rank of $\La$:
\[
\operatorname{rank}(\hat{L}) = \operatorname{rank}(\La).
\]

\item[(3)] The concentration matrix $\hat{S}-\hat{L}$ can be realized as the
marginal concentration matrix of an appropriate latent-variable model:
\[
\hat{S} - \hat{L} \succ0,\qquad \hat{L} \succeq0.
\]
\end{longlist}
\end{defin}

When a sequence of estimators is algebraically correct with probability
approaching 1 in a suitable high-dimensional scaling regime, then we
say\vspace*{1pt} that the estimators are \textit{algebraically
consistent}. The first condition ensures that $\hat{S}$ provides the
correct structural estimate of the conditional graphical model of the
observed variables conditioned on the latent components. This property
is the same as the ``sparsistency'' property studied in standard
graphical model selection~\cite{LamF2009,RavWRY2008}. The second
condition ensures that the number of latent\vspace*{1pt} components is properly
estimated. Finally, the third condition ensures that the pair of
matrices $(\hat{S},\hat{L})$ leads to a realizable latent-variable
model. In particular, this condition implies that there exists a valid
latent-variable model in which (a) the conditional graphical model
structure among the observed variables is given by $\hat{S}$, \textup{(b)} the
number of latent variables is equal to the rank of $\hat{L}$, and (c)
the extra correlations induced due to marginalization over the latent
variables are equal to $\hat{L}$. Any method for matrix factorization
(e.g.,~\cite{WitTH2009}) can be used to further factorize $\hat{L}$,
depending on the property that one desires in the factors (e.g.,
sparsity).

We also study estimation error rates in the usual sense, that is, we
show that one can produce estimates $(\hat{S},\hat{L})$ that are close
in various norms to the matrices $(\Sa,\La)$. Notice that bounding the
estimation error in some norm does not in general imply that the
support/sign-pattern and rank of $(\hat{S},\hat{L})$ are the same as
those of $(\Sa,\La)$. Therefore bounded estimation error is different
from algebraic correctness, which requires that $(\hat{S},\hat{L})$
have the same support/sign-pattern and rank as $(\Sa,\La)$.

\subsubsection*{Goal} Let $K^\ast_{(\OH)}$ denote the concentration
matrix of a Gaussian model. Suppose that we have $n$ samples $\{X^i_O\}
_{i=1}^n$\vspace*{1pt} of the observed variables $X_O$. We would like to produce
estimates $(\hat{S}_n,\hat{L}_n)$ that, with high probability, are
algebraically correct and have bounded estimation error (in some norm).

\subsubsection*{Our approach} We propose\vspace*{1pt} the regularized likelihood
convex program (\ref{eqsdp}) to produce estimates $(\hat{S}_n,\hat
{L}_n)$. Specifically, the sample covariance matrix $\Sigma_O^n$ in
(\ref{eqsdp}) is defined as
\[
\Sigma^n_O \triangleq\frac{1}{n}\sum_{i=1}^n X^i_O {X^i_O}^T.
\]
We give conditions on the underlying model $K^\ast_{(\OH)}$ and
suitable choices for the parameters $\lambda_n, \gamma$ under which the
estimates $(\hat{S}_n,\hat{L}_n)$ are consistent (see Theorem~\ref
{theomain}).

\subsection{Likelihood function and Fisher information}
\label{subsecll}


$\!\!$Given $n$ samples $\{X^i\}_{i=1}^n$ of a finite collection of jointly
Gaussian zero-mean random variables with concentration matrix $K^\ast$,
it is easily seen that the log-likelihood function is given by:
%
\begin{equation}
\lh(K; \Sigma^n) = \log\det(K) - \tr(K \Sigma^n),
\end{equation}
where $\lh(K;\Sigma^n)$ is a function of $K$. Notice that this function
is strictly concave for $K \succ0$. In the latent-variable modeling
problem with sample covariance $\Sigma_O^n$, the likelihood function
with respect to the parametrization $(S,L)$ is given by $\lh(S-L;
\Sigma_O^n)$. This function is \textit{jointly concave} with respect to
the parameters $(S,L)$ whenever $S-L \succ0$, and it is employed in
our variational formulation~(\ref{eqsdp}) to learn a latent-variable model.

In the analysis of a convex program involving the likelihood function,
the Fisher information plays an important role as it is the negative of
the Hessian of the likelihood function and thus controls the curvature.
As the first term in the likelihood function is linear, we need only
study higher-order derivatives of the log-determinant function in order
to compute the Hessian. In the latent-variable setting with the
marginal concentration matrix of the observed variables given by
$\tilde
{K}^\ast_O = (\Sigma^\ast_O)^{-1}$ [see~(\ref{eqschur})], the
corresponding Fisher information matrix is
%
\begin{equation} \label{eqfisher}
\F(\tilde{K}^\ast_O) = (\tilde{K}^\ast_O)^{-1} \otimes(\tilde
{K}^\ast
_O)^{-1} = \Sigma^\ast_O \otimes\Sigma^\ast_O.
\end{equation}
Here $\otimes$ denotes the tensor product between matrices. Notice that
this is precisely the $|O|^2 \times|O|^2$ submatrix of the full Fisher
information matrix $\F(K^\ast_{(\OH)}) = \Sigma_{(\OH)}^\ast\otimes
\Sigma_{(\OH)}^\ast$ with respect to all the parameters $K^\ast_{(\OH)}
= (\Sigma_{(\OH)}^\ast)^{-1}$ (corresponding to the situation in which
\textit{all} the variables $X_{O \cup H}$ are observed). In Section \ref
{subsecfi} we impose various conditions on the Fisher information
matrix $\F(\tilde{K}^\ast_O)$ under which our regularized
maximum-likelihood formulation provides consistent estimates.

\subsection{Algebraic varieties of sparse and low-rank matrices}
\label{subsecav}

The set of sparse matrices and the set of low-rank matrices can be
naturally viewed as algebraic varieties (solution sets of systems of
polynomial equations). Here we describe these varieties, and discuss
some of their geometric properties such as the tangent space and local
curvature at a (smooth) point.

Let $\s(k)$ denote the set of matrices with at most $k$ nonzeros:
%
\begin{equation}
\s(k) \triangleq\{M \in\mathbb{R}^{p \times p}  |  |
{\operatorname{support}}(M)| \leq k \}.
\end{equation}
Here $\mathrm{support}$ denotes the locations of nonzero entries. The
set $\s(k)$ is an algebraic variety, and can in fact be viewed as a
union of ${p^2 \choose k}$\vadjust{\goodbreak} subspaces in $\mathbb{R}^{p \times p}$. This
variety has dimension $k$, and it is smooth everywhere except at those
matrices that have support size strictly smaller than $k$. For any
matrix $M \in\mathbb{R}^{p \times p}$, consider the variety $\s
(|\mathrm{support}(M)|)$; $M$ is a smooth point of this variety, and
the tangent space at $M$ is given by
%
\begin{equation}
\Omega(M) = \{N \in\mathbb{R}^{p \times p}  |  \mathrm{support}(N)
\subseteq\operatorname{support}(M) \}.
\end{equation}


Next let $\lr(r)$ denote the algebraic variety of matrices with rank at
most~$r$:
%
\begin{equation}
\lr(r) \triangleq\{M \in\mathbb{R}^{p \times p}  |
{\operatorname{rank}}(M) \leq r \}.
\end{equation}
It is easily seen that $\lr(r)$ is an algebraic variety because it can
be defined through the vanishing of all $(r+1) \times(r+1)$ minors.
This variety has dimension equal to $r (2p - r)$, and it is smooth
everywhere except at those matrices that have rank strictly smaller
than $r$. Consider a rank-$r$ matrix $M$ with singular value
decomposition (SVD) given by $M = U D V^T$, where $U,V \in\mathbb
{R}^{p \times r}$ and $D \in\mathbb{R}^{r \times r}$. The matrix $M$
is a smooth point of the variety $\lr(\operatorname{rank}(M))$, and the
tangent space at $M$ with respect to this variety is given by
%
\begin{equation}
T(M) = \{U Y_1^T + Y_2 V^T  |  Y_1,Y_2 \in\mathbb{R}^{p \times r}\}.
\end{equation}
We view both $\Omega(M)$ and $T(M)$ as subspaces in $\mathbb{R}^{p
\times p}$. In Section~\ref{seciden} we explore the connection between
geometric properties of these tangent spaces and the identifiability
problem in latent-variable graphical models.


\subsubsection*{Curvature of rank variety} The sparse matrix variety
$\s
(k)$ has the property that it has \textit{zero} curvature at any smooth
point. The situation is more complicated for the low-rank matrix
variety $\lr(r)$, because the curvature at any smooth point is nonzero.
We analyze how this variety curves locally, by studying how the tangent
space changes from one point to a neighboring point. Indeed the amount
of curvature at a point is directly related to the ``angle'' between
the tangent space at that point and the tangent space at a neighboring
point. For any linear subspace $T$ of matrices, let $\proj_{T}$ denote
the projection onto $T$. Given two subspaces $T_1,T_2$ of the same
dimension, we measure the ``twisting'' between these subspaces by
considering the following quantity:
%
\begin{equation}\label{eqrho}
\rho(T_1,T_2) \triangleq\|\proj_{T_1} - \proj_{T_2}\|_{2
\rightarrow
2} = {\max_{\|N\|_2 \leq1} } \|[\proj_{T_1} - \proj_{T_2}] (N)\|_2.
\end{equation}
In the supplement~\cite{supp} we review relevant results from matrix
perturbation theory, which suggest that the magnitude of the smallest
nonzero singular value is closely tied to the local curvature of the
variety. Therefore we control the twisting between tangent spaces at
nearby points by bounding the smallest nonzero singular value away from zero.


%
%
%
%
%


\section{Identifiability}
\label{seciden}

In the absence of additional conditions, the latent-variable model
selection problem is ill-posed. In this section we discuss a set of
conditions on latent-variable models that ensure\vadjust{\goodbreak} that these models are
identifiable given marginal statistics for a subset of the variables.
Some of the discussion in Sections~\ref{subseclv} and \ref
{subsecov} is presented in greater detail
in~\cite{ChaSPW2009}.\looseness=1

\subsection{Structure between latent and observed variables}
\label{subseclv}

Suppose that the low-rank matrix that summarizes the effect of the
latent components is itself sparse. This leads to identifiability
issues in the sparse-plus-low-rank decomposition problem. Statistically
the additional correlations induced due to marginalization over the
latent variables could be mistaken for the conditional graphical model
structure of the observed variables. In order to avoid such
identifiability problems the effect of the latent variables must be
``diffuse'' across the observed variables. To address this point the
following quantity was introduced in~\cite{ChaSPW2009} for any matrix
$M$, defined with respect to the tangent space $T(M)$:
%
\begin{equation} \label{eqxi}
\xi(T(M)) \triangleq{\max_{N \in T(M),  \|N\|_2 \leq1}}
\|N\|_\infty.
\end{equation}
Thus $\xi(T(M))$ being small implies that elements of the tangent
space $T(M)$ cannot have their support concentrated
in a few locations; as a result $M$ cannot be too sparse. This idea is
formalized in~\cite{ChaSPW2009} by relating $\xi(T(M))$ to a notion of
``incoherence'' of the row/column spaces, where the row/column spaces
are said to be incoherent with respect to the standard basis if these
spaces are not aligned closely with any of the coordinate axes.
Typically a matrix $M$ with
incoherent row/column spaces would have $\xi(T(M)) \ll1$. This point
is quantified precisely in~\cite{ChaSPW2009}. Specifically, we note
that $\xi(T(M))$ can be as small as $\sim\sqrt{\frac{r}{p}}$ for a
rank-$r$ matrix $M \in\mathbb{R}^{p \times p}$ with row/column spaces
that are almost maximally incoherent (e.g., if the row/column spaces
span any $r$ columns of a $p \times p$ orthonormal Hadamard matrix). On
the other hand, $\xi(T(M)) = 1$ if the row/column spaces of $M$ contain
a standard basis vector.

Based on these concepts we roughly require that the low-rank matrix
that summarizes the effect of the latent variables be \textit
{incoherent}, thereby ensuring that the extra correlations due to
marginalization over the latent components cannot be confused with the
conditional graphical model structure of the observed variables. Notice
that the quantity $\xi$ is not just a measure of the number of latent
variables, but also of the overall effect of the correlations induced
by marginalization over these variables.

\textit{Curvature and change in $\xi$}: As noted previously, an
important technical point is that the algebraic variety of low-rank
matrices is locally curved at any smooth point. Consequently the
quantity $\xi$ changes as we move along the low-rank matrix variety
smoothly. The quantity $\rho(T_1,T_2)$ introduced in~(\ref{eqrho})
allows us to bound the variation in $\xi$ as follows (proof in
Section~\ref{secproofs}):

\begin{lemm}\label{theorhotspace}
Let $T_1,T_2$ be two linear subspaces of matrices of the same dimension
with the property that $\rho(T_1,T_2) < 1$, where $\rho$ is defined in
(\ref{eqrho}). Then we have that
\[
\xi(T_2) \leq\frac{1}{1-\rho(T_1,T_2)}  [\xi(T_1) + \rho(T_1,T_2)].
\]
\end{lemm}

\subsection{Structure among observed variables}
\label{subsecov}

An identifiability problem also arises if the conditional graphical
model among the observed variables contains a densely connected
subgraph. These statistical relationships might be mistaken as
correlations induced by marginalization over latent variables.
Therefore we need to ensure that the conditional graphical model among
the observed variables is sparse. We impose the condition that this
conditional graphical model must have small ``degree,'' that is, no
observed variable is directly connected to too many other observed
variables conditioned on the latent components. Notice that bounding
the degree is a more refined condition than simply bounding the total
number of nonzeros as the \textit{sparsity pattern} also plays a role. In
\cite{ChaSPW2009} the authors introduced the following quantity in
order to provide an appropriate measure of the sparsity pattern of a matrix:
%
\begin{equation} \label{eqmu}
\mu(\Omega(M)) \triangleq{\max_{N \in\Omega(M), \|N\|_\infty\leq1}}
\|N\|_2.
\end{equation}
The quantity $\mu(\Omega(M))$ being small for a matrix implies that
the spectrum of any element of the tangent space $\Omega(M)$
is not too ``concentrated,'' that is, the singular values of the
elements of the tangent space are not too large. In~\cite{ChaSPW2009}
it is shown that a sparse matrix $M$ with ``bounded degree'' (a small
number of nonzeros per row/column) has small $\mu(M)$. Specifically, if
$M \in\mathbb{R}^{p \times p}$ is any matrix with at most $\degr(M)$
nonzero entries per row/column, then we have that
\[
\mu(\Omega(M)) \leq\degr(M).
\]

\subsection{Transversality of tangent spaces}
\label{subsectts}

Suppose that we have the sum of two vectors, each from two known
subspaces. It is possible to uniquely recover the individual vectors
from the sum if and only if the subspaces have a transverse
intersection, that is, they only intersect at the origin. This simple
observation leads to an appealing geometric notion of identifiability.
Suppose now that we have the sum of a sparse matrix and a low-rank
matrix, and that we are also given the tangent spaces at these matrices
with respect to the algebraic varieties of sparse and low-rank
matrices, respectively. Then a necessary and sufficient condition for
identifiability with respect to the tangent spaces is that these spaces
have a transverse intersection. This transverse intersection condition
is also sufficient for local identifiability in a neighborhood around
the sparse matrix and low-rank matrix with respect to the varieties of
sparse and low-rank matrices (due to the inverse function theorem). It
turns out that these tangent space transversality conditions are also
sufficient for the convex program~(\ref{eqsdp}) to provide consistent
estimates of a latent-variable graphical model (without any side
information about the tangent spaces).

In order to quantify the level of transversality between the tangent
spaces $\Omega$ and $T$ we study the \textit{minimum gain} with respect
to some norm of the addition operator (which adds two matrices)
$\add\dvtx
\mathbb{R}^{p \times p} \times\mathbb{R}^{p \times p} \rightarrow
\mathbb{R}^{p \times p}$ restricted to the cartesian product $\ts=
\Omega\times T$. Then given any matrix norm \mbox{$\|\cdot\|_q$} on $\mathbb
{R}^{p \times p} \times\mathbb{R}^{p \times p}$, the minimum gain of
$\add$ restricted to $\ts$ is defined as
\[
\varepsilon(\Omega,T,\|\cdot\|_q) \triangleq{\min_{(S,L) \in\Omega
\times
T,  \|(S,L)\|_q = 1} } \|\proj_\ts\add^\dag\add\proj_\ts
(S,L)\|_q,
\]
where $\proj_\ts$ denotes the projection onto $\ts$, and $\add^\dag$
denotes the adjoint of the addition operator (with respect to the
standard Euclidean inner-product). The ``level'' of transversality of
$\Omega$ and $T$ is measured by the magnitude of $\varepsilon(\Omega
,T,\|
\cdot\|_q)$, with transverse intersection being equivalent to
$\varepsilon
(\Omega,T,\|\cdot\|_q) > 0$. Note that $\varepsilon(\Omega,T,\|\cdot\|_F)$
is the square of the \textit{minimum singular value} of the addition
operator $\add$ restricted to $\Omega\times T$.

A natural norm with which to measure transversality is the dual norm of
the regularization function in~(\ref{eqsdp}), as the subdifferential
of the regularization function is specified in terms of its dual. The
reasons for this will become clearer as we proceed through this paper.
Recall that the regularization function used in the variational
formulation~(\ref{eqsdp}) is given by
\[
f_\gamma(S,L) = \gamma\|S\|_1 + \|L\|_\ast,
\]
where the nuclear norm \mbox{$\|\cdot\|_\ast$} reduces to the trace function
over the cone of positive-semidefinite matrices. This function is a
norm for all $\gamma> 0$. The dual norm of $f_\gamma$ is given by
\[
\g(S,L) = \max\biggl\{\frac{\|S\|_\infty}{\gamma}, \|L\|_2
\biggr\}.
\]

%

Next we define the quantity $\chi(\Omega,T,\gamma)$ as follows in order
to study the transversality of the spaces $\Omega$ and $T$ with respect
to the $\g$ norm:
%
\begin{equation} \label{eqchi}
\chi(\Omega,T,\gamma) \triangleq\max\biggl\{\frac{\xi(T)}{\gamma
}, 2 \mu
(\Omega) \gamma\biggr\}.
\end{equation}
Here $\mu$ and $\xi$ are defined in~(\ref{eqmu}) and (\ref
{eqxi}). We
then have the following result (proved in Section~\ref{secproofs}):
%
\begin{lemm}\label{theorg1}
Let $S \in\Omega, L \in T$ be matrices such that $\|S\|_{\infty} =
\gamma$ and let $\|L\|_2 = 1$. Then we have that $\g(\proj_{\ts}
\add
^\dag\add\proj_\ts(S,L)) \in[1-\chi(\Omega,T,\gamma),1+\chi
(\Omega
,T,\gamma)]$, where $\ts= \Omega\times T$ and $\chi(\Omega
,T,\gamma)$
is defined in~(\ref{eqchi}). In particular we have that $1-\chi
(\Omega
,T,\gamma) \leq\varepsilon(\Omega,T,\g)$.
\end{lemm}

The quantity $\chi(\Omega,T,\gamma)$ being small implies that the
addition operator is essentially isometric when restricted to $\ts=
\Omega\times T$. Stated differently, the magnitude of $\chi(\Omega
,T,\gamma)$ is a measure of the level of transversality of the spaces
$\Omega$ and $T$. If $\mu(\Omega) \xi(T) < \frac{1}{2}$, then
$\gamma
\in(\xi(T), \frac{1}{2 \mu(\Omega)})$ ensures\vspace*{1pt} that $\chi(\Omega
,T,\gamma) < 1$, which in turn implies that the tangent spaces $\Omega$
and $T$ have a transverse intersection.

\textit{Observation}: Thus we have that the smaller the quantities
$\mu
(\Omega)$ and $\xi(T)$, the more transverse the intersection of the
spaces $\Omega$ and $T$ as measured by $\varepsilon(\Omega,T,\g)$.

\subsection{Conditions on Fisher information}
\label{subsecfi}

The main focus of Section~\ref{secmain} is to analyze the regularized
maximum-likelihood convex program~(\ref{eqsdp}) by studying its
optimality conditions. The log-likelihood function is well-approximated
in a neighborhood by a quadratic form given by the Fisher information
(which measures the curvature, as discussed in Section \ref
{subsecll}). Let $\F^\ast= \F(\tilde{K}^\ast_O)$ denote the Fisher
information evaluated at the true marginal concentration matrix $\tilde
{K}^\ast_O$ [see~(\ref{eqschur})]. The appropriate measure of
transversality between the tangent
spaces\setcounter{footnote}{1}\footnote{We implicitly assume
that these tangent spaces are subspaces of the space of \textit
{symmetric} matrices.} $\Omega= \Omega(\Sa)$ and $T = T(\La)$ is then
in a space in which the inner-product is given by $\F^\ast$.
Specifically, we need to analyze the minimum gain of the operator
$\proj
_{\ts} \add^\dag\F^\ast\add\proj_{\ts}$ restricted to the space
$\ts
= \Omega\times T$. Therefore we impose several conditions on the
Fisher information $\F^\ast$. We define quantities that control the
gains of $\F^\ast$ restricted to $\Omega$ and $T$ separately; these
ensure that elements of $\Omega$ and elements of $T$ are individually
identifiable under the map $\F^\ast$. In addition we define quantities
that, in conjunction with bounds on $\mu(\Omega)$ and $\xi(T)$, allow
us to control the gain of $\F^\ast$ restricted to the direct-sum
$\Omega\oplus T$.

\textit{$\F^\ast$ restricted to $\Omega$}: The minimum gain of the
operator $\proj_{\Omega} \F^\ast\proj_{\Omega}$ restricted to
$\Omega$
is given by
\[
\alpha_{\Omega} \triangleq{\min_{M \in\Omega, \|M\|_\infty= 1}}
\|
\proj_{\Omega} \F^\ast\proj_{\Omega}(M)\|_\infty.
\]
The maximum effect of elements in $\Omega$ in the orthogonal direction
$\Omega^\bot$ is given by
\[
\delta_{\Omega} \triangleq{\max_{M \in\Omega, \|M\|_\infty= 1}}
\|
\proj_{\Omega^\bot} \F^\ast\proj_{\Omega}(M)\|_\infty.
\]
The operator $\F^\ast$ is injective on $\Omega$ if $\alpha_\Omega> 0$.
The ratio $\frac{\delta_\Omega}{\alpha_\Omega} \leq1 - \nu$ implies
the irrepresentability condition imposed in~\cite{RavWRY2008}, which
gives a sufficient condition for consistent recovery of graphical model
structure using $\ell_1$-regularized maximum-likelihood. Notice that
this condition is a generalization of the usual Lasso
irrepresentability conditions~\cite{ZhaY2006,Wai2009}, which are
typically imposed on the covariance matrix. Finally we also consider
the following quantity, which controls the behavior of $\F^\ast$
restricted to $\Omega$ in the spectral norm:
\[
\beta_\Omega\triangleq{\max_{M \in\Omega, \|M\|_2 = 1}}
\|\F^\ast(M)\|_2.
\]

\textit{$\F^\ast$ restricted to $T$}: Analogously to the case of
$\Omega
$ one could control the gains of the operators $\proj_{T^\bot} \F
^\ast
\proj_T$ and $\proj_T \F^\ast\proj_T$. However, as discussed
previously, one complication is that the tangent spaces at nearby
smooth points on the rank variety are in general different, and the
amount of twisting between these spaces is governed by the local
curvature. Therefore we control the gains of the operators $\proj
_{T'^\bot} \F^\ast\proj_{T'}$ and $\proj_{T'} \F^\ast\proj
_{T'}$ for
all tangent spaces $T'$ that are ``close to'' the nominal $T$ (at the
true underlying low-rank matrix), measured by $\rho(T,T')$ (\ref
{eqrho}) being small. The minimum gain of the operator $\proj_{T'} \F
^\ast\proj_{T'}$ restricted to $T'$ (close to $T$) is given by
\[
\alpha_{T} \triangleq{\min_{\rho(T',T) \leq{\xi(T)}/{2}}
\min_{M
\in T', \|M\|_2 = 1}} \|\proj_{T'} \F^\ast\proj_{T'}(M)\|_2.
\]
Similarly, the maximum effect of elements in $T'$ in the orthogonal
direction $T'^\bot$ (for $T'$ close to $T$) is given by
\[
\delta_{T} \triangleq{\max_{\rho(T',T) \leq{\xi(T)}/{2}}
\max_{M
\in T', \|M\|_2 = 1} } \|\proj_{T'^\bot} \F^\ast\proj_{T'}(M)\|_2.
\]
Implicit in the definition of $\alpha_T$ and $\delta_T$ is the fact
that the outer minimum and maximum are only taken over spaces $T'$ that
are tangent spaces to the rank-variety. The operator $\F^\ast$ is
injective on all tangent spaces $T'$ such that $\rho(T',T) \leq\frac
{\xi(T)}{2}$ if $\alpha_T > 0$. An irrepresentability condition
(analogous\vspace*{1pt} to those developed for the sparse case) for tangent spaces
near $T$ to the rank variety would be that $\frac{\delta_T}{\alpha_T}
\leq1 - \nu$. Finally we also control the behavior of $\F^\ast$
restricted to $T'$ close to $T$ in the $\ell_\infty$ norm:
\[
\beta_T \triangleq{\max_{\rho(T',T) \leq{\xi(T)}/{2}}  \max_{M
\in T', \|M\|_\infty= 1}}  \|\F^\ast(M)\|_\infty.
\]

The two sets of quantities $(\alpha_\Omega,\delta_\Omega)$ and
$(\alpha
_T,\delta_T)$ essentially control how $\F^\ast$ behaves when restricted
to the spaces $\Omega$ and $T$ \textit{separately} (in the natural
norms). The quantities $\beta_\Omega$ and $\beta_T$ are useful in order
to control the gains of the operator $\F^\ast$ restricted to the
\textit
{direct sum} $\Omega\oplus T$. Notice that although the magnitudes of
elements in $\Omega$ are measured most naturally in the $\ell_\infty$
norm, the quantity $\beta_\Omega$ is specified with respect to the
spectral norm. Similarly, elements of the tangent spaces $T'$ to the
rank variety are most naturally measured in the spectral norm, but
$\beta_T$ provides control in the $\ell_\infty$ norm. These quantities,
combined with $\mu(\Omega)$ and $\xi(T)$ [defined in~(\ref{eqmu}) and
(\ref{eqxi})], provide the ``coupling'' necessary to control the
behavior of $\F^\ast$ restricted to elements in the direct sum
$\Omega
\oplus T$. In order to keep track of fewer quantities, we summarize the
six quantities as follows:
\[
\alpha\triangleq\min(\alpha_\Omega, \alpha_T); \qquad \delta
\triangleq
\max(\delta_\Omega, \delta_T); \qquad \beta\triangleq\max(\beta
_\Omega,
\beta_T).
\]

\textit{Main assumption}: There exists a $\nu\in(0, \frac
{1}{2}]$ such that
\[
\frac{\delta}{\alpha} \leq1 - 2 \nu.
\]

This assumption is to be viewed as a generalization of the
irrepresentability conditions imposed on the covariance matrix \cite
{ZhaY2006,Wai2009} or the Fisher information matrix~\cite{RavWRY2008}
in order to provide consistency guarantees for sparse model selection
using the $\ell_1$ norm. With this assumption we have the following
proposition, proved in Section~\ref{secproofs}, about the gains of the
operator $\F^\ast$ restricted to $\Omega\oplus T$. This proposition
plays a fundamental role in the analysis of the performance of the
regularized maximum-likelihood procedure~(\ref{eqsdp}). Specifically,
it gives conditions under which a suitable primal-dual pair can be
specified to certify optimality with respect to~(\ref{eqsdp}) (see
Section~\ref{subsecproof} for more details).
%
\begin{prop}\label{theoirr}
Let $\Omega$ and $T$ be the tangent spaces defined in this section, and
let $\F^\ast$ be the Fisher information evaluated at the true marginal
concentration matrix. Further let $\alpha,\beta,\nu$ be as defined
above. Suppose that
\[
\mu(\Omega) \xi(T) \leq\frac{1}{6} \biggl(\frac{\nu\alpha
}{\beta(2 -
\nu)}\biggr)^2,
\]
and that $\gamma$ is in the following range:
\[
\gamma\in\biggl[\frac{3 \xi(T) \beta(2-\nu)}{\nu\alpha}, \frac
{\nu
\alpha}{2 \mu(\Omega) \beta(2-\nu)} \biggr].
\]
Then we have the following two conclusions for $\ts= \Omega\times T'$
with $\rho(T', T) \leq\frac{\xi(T)}{2}$:
\begin{longlist}[(2)]
\item[(1)] The minimum gain of $\F^\ast$ restricted to $\Omega\oplus T'$ is
bounded below:
\[
\min_{(S,L) \in\ts,  \|S\|_\infty= \gamma,  \|L\|_2 = 1}  \g
(\proj
_\ts\add^\dag\F^\ast\add\proj_\ts(S,L)) \geq \frac{\alpha}{2}.
\]
Specifically this implies that for all $(S,L) \in\ts$
\[
\g(\proj_\ts\add^\dag\F^\ast\add\proj_\ts(S,L)) \geq\frac
{\alpha
}{2} \g(S,L).
\]

\item[(2)] The effect of elements in $\ts= \Omega\times T'$ on the
orthogonal complement $\ts^\bot= \Omega^\bot\times T'^\bot$ is
bounded above:
\[
\|\proj_{\ts^\bot} \add^\dag\F^\ast\add\proj_\ts
(\proj
_\ts\add^\dag\F^\ast\add\proj_\ts)^{-1} \|_{\g
\rightarrow\g} \leq1-\nu.
\]
Specifically this implies that for all $(S,L) \in\ts$
\[
\g(\proj_{\ts^\bot} \add^\dag\F^\ast\add\proj_\ts(S,L)) \leq
(1-\nu
) \g(\proj_\ts\add^\dag\F^\ast\add\proj_\ts(S,L)).
\]
\end{longlist}
\end{prop}



The last quantity we consider is the spectral norm of the marginal
covariance matrix $\Sigma^\ast_O = (\tilde{K}^\ast_O)^{-1}$:
%
\begin{equation}\label{eqpsi}
\psi\triangleq\|\Sigma^\ast_O\|_2 = \|(\tilde{K}^\ast_O)^{-1}\|_2.
\end{equation}
A bound on $\psi$ is useful in the probabilistic component of our
analysis, in order to derive convergence rates of the sample covariance
matrix to the true covariance matrix. We also observe that
\[
\|\F^\ast\|_{2 \rightarrow2} = \|(\tilde{K}^\ast_O)^{-1} \otimes
(\tilde{K}^\ast_O)^{-1}\|_{2 \rightarrow2} = \psi^2.
\]

\subsubsection*{Remarks} The quantities $\alpha,\beta,\delta$ bound the
gains of the Fisher information $\F^\ast$ restricted to the spaces
$\Omega$ and $T$ (and tangent spaces near $T$). One can make stronger
assumptions on $\F^\ast$ that are more easily interpretable. For
example, $\alpha_\Omega, \beta_\Omega$ could bound the minimum/maximum
gains of $\F^\ast$ for \textit{all} matrices (rather than just those in
$\Omega$), and $\delta_\Omega$ the $\F^\ast$-inner-product for
\textit
{all} pairs of orthogonal matrices (rather than just those in $\Omega$
and $\Omega^\bot$). Similarly, $\alpha_T,\beta_T$ could bound the
minimum/maximum gains of $\F^\ast$ for all matrices (rather than just
those near $T$), and $\delta_T$ the $\F^\ast$-inner-product for all
pairs of orthogonal matrices (rather than just those near $T$ and
$T^\bot$). Such bounds would apply in either the $\|\cdot\|_{2
\rightarrow2}$ norm (for $\alpha_T, \delta_T, \beta_\Omega$) or
the $\|
\cdot\|_{\infty\rightarrow\infty}$ norm (for $\alpha_\Omega,
\delta
_\Omega, \beta_T$). These modified assumptions are global in nature
(not restricted just to $\Omega$ or near $T$) and are consequently
stronger (they lower-bound the original $\alpha_\Omega,\alpha_T$ and
they upper-bound the original $\beta_\Omega,\beta_T,\delta_\Omega
,\delta
_T$), and they essentially control the gains of the operator $\F^\ast$
in the $\|\cdot\|_{2 \rightarrow2}$ norm and the $\|\cdot\|_{\infty
\rightarrow\infty}$ norm. In contrast, previous works on covariance
selection~\cite{BicL2008a,BicL2008b,RotBLZ2008} consider \textit
{well-conditioned} families of covariance matrices by bounding the
minimum/maximum eigenvalues (i.e., gain with respect to the spectral norm).


\section{Consistency of regularized maximum-likelihood program}
\label{secmain}

\subsection{Main results}
\label{subsecmainres}

Recall\vspace*{1pt} that $K^\ast_{(\OH)}$ denotes the full concentration matrix of a
collection of zero-mean jointly-Gaussian observed and latent variables.
Let $p = |O|$ denote the number of observed variables, and let $h =
|H|$ denote the number of latent variables. We are given $n$ samples $\{
X_O^i\}_{i=1}^n$ of the observed variables $X_O$. We consider the
high-dimensional setting in which $(p,h,n)$ are all allowed to grow
simultaneously. We present our main result next demonstrating the
consistency of the estimator~(\ref{eqsdp}), and then discuss classes
of latent-variable graphical models and various scaling regimes in
which our estimator is consistent. Recall from~(\ref{eqsdp}) that
$\lambda_n$ is a regularization parameter, and $\gamma$ is a trade-off
parameter between the rank and sparsity terms. Notice from
Proposition~\ref{theoirr} that the choice of $\gamma$ depends on the
values of $\mu(\Omega(\Sa))$ and $\xi(T(\La))$. While these quantities
may not be known a priori, we discuss a method to choose $\gamma$
numerically in our experimental results (see Section~\ref{secsims}).
The following theorem shows that the estimates $(\hat{S}_n, \hat{L}_n)$
provided by the convex program~(\ref{eqsdp}) are consistent for a
suitable choice of $\lambda_n$. In addition to the appropriate
identifiability conditions (as specified by Proposition \ref
{theoirr}), we also impose lower bounds on the minimum magnitude
nonzero entry $\theta$ of the sparse conditional graphical model matrix
$\Sa$ and on the minimum nonzero singular value $\sigma$ of the
low-rank matrix $\La$ summarizing the effect of the latent variables.
The theorem is stated in terms of the quantities $\alpha,\beta,\nu
,\psi
$, and we particularly emphasize the dependence on $\mu(\Omega(\Sa))$
and $\xi(T(\La))$ because these control the complexity of the
underlying latent-variable graphical model given by $K^\ast_{(\OH)}$. A
number of quantities play a role in our theorem: let\vspace*{-1pt} $D = \max\{
1,\frac
{\nu\alpha}{3 \beta(2-\nu)}\}$, $C_1 = \psi(1 + \frac{\alpha}{6
\beta})$, $C_2 = \frac{48}{\alpha} + \frac{1}{\psi^2}$,
$C_{\mathrm
{samp}} = \frac{\alpha\nu}{32 (3-\nu) D} \min\{ \frac{1}{4
C_1},\frac{\alpha\nu}{256 D (3-\nu) \psi C_1^2} \}$,\vspace*{-1pt}
$C_\lambda
= \frac{48 \sqrt{2} D \psi(2-\nu)}{\xi(T) \nu}$, $C_S = \max
\{
(\frac{6(2-\nu)}{\nu} + 1 ) C_2^2 \psi^2 D, C_2 +
\frac{3
\alpha C_2^2 (2-\nu)}{16 (3-\nu)} \}$ and $C_L = \frac{C_2
\nu
\alpha}{\beta(2-\nu)}$.
%
%
%
\begin{theorem} \label{theomain}
Let $K^\ast_{(\OH)}$ denote the concentration matrix of a Gaussian
model. We have $n$ samples $\{X^i_O\}_{i=1}^n$ of the $p$ observed
variables denoted by $O$. Let $\Omega= \Omega(\Sa)$ and $T = T(\La)$
denote the tangent spaces at $\Sa$ and at $\La$ with respect to the
sparse and low-rank matrix varieties, respectively.


\textup{Assumptions}: Suppose that the quantities $\mu(\Omega)$ and
$\xi
(T)$ satisfy the assumption of Proposition~\ref{theoirr} for
identifiability, and $\gamma$ is chosen in the range specified by
Proposition~\ref{theoirr}. Further suppose that the following
conditions hold:

\begin{longlist}[(4)]
\item[(1)] Let $n \geq\frac{p}{\xi(T)^4} \max\{\frac{128 \psi
^2}{C^2_{\mathrm{samp}}}, 2\}$, that is, we require that $n
\gtrsim\frac{p}{\xi(T)^4}$.

\item[(2)] Set $\lambda_n = \frac{48 \sqrt{2} D \psi(2-\nu)}{\xi(T)
\nu}
\sqrt{\frac{p}{n}}$, that is, we require that $\lambda_n \asymp
\frac
{1}{\xi(T)} \sqrt{\frac{p}{n}}$.

\item[(3)] Let $\sigma\geq\frac{C_L \lambda_n}{\xi(T)^2}$, that is, we
require that $\sigma\gtrsim\frac{1}{\xi(T)^3} \sqrt{\frac{p}{n}}$.

\item[(4)] Let $\theta\geq\frac{C_S \lambda_n}{\mu(\Omega)}$, that
is, we
require that $\theta\gtrsim\frac{1}{\xi(T) \mu(\Omega)} \sqrt
{\frac{p}{n}}$.
\end{longlist}

\textup{Conclusions}: Then with probability greater than $1 - 2\exp\{
-p\}$ we have algebraic correctness and estimation error given
by:
\begin{longlist}[(2)]
\item[(1)] $\operatorname{sign}(\hat{S}_n) = \operatorname{sign}(\Sa)$ and
$\operatorname{rank}(\hat{L}_n) = \operatorname{rank}(\La)$;

\item[(2)] $\g(\hat{S}_n-\Sa, \hat{L}_n - \La) \leq\frac{512 \sqrt{2}
(3-\nu) D \psi}{\nu\alpha\xi(T)} \sqrt{\frac{p}{n}} \lesssim
\frac
{1}{\xi(T)} \sqrt{\frac{p}{n}}$.
\end{longlist}
\end{theorem}

The proof of this theorem is given in Section~\ref{secproofs}. The
theorem essentially states that if the minimum nonzero singular value
of the low-rank piece $\La$ and minimum nonzero entry of the sparse
piece $\Sa$ are bounded away from zero, then the convex program (\ref
{eqsdp}) provides estimates that are both algebraically correct and
have bounded estimation error (in the $\ell_\infty$ and spectral norms).


Notice that the condition on the minimum singular value of $\La$ is
more stringent than the one on the minimum nonzero entry of $\Sa$. One
role played by these conditions is to ensure that the estimates $(\hat
{S}_n,\hat{L}_n)$ do not have smaller support size/rank than $(\Sa
,\La
)$. However, the minimum singular value bound plays the additional role
of bounding the curvature of the low-rank matrix variety around the
point $\La$, which is the reason for this condition being more
stringent. Notice also that the number of latent variables $h$ does not
explicitly appear in the bounds in Theorem~\ref{theomain}, which only
depend on $p,\mu(\Omega(\Sa)),\xi(T(\La))$. However, the
dependence on
$h$ is implicit in the dependence on $\xi(T(\La))$, and we discuss this
point in greater detail in the following section.

Finally we note that consistency holds in Theorem~\ref{theomain} for a
\textit{range} of values of $\gamma\in[\frac{3 \beta(2-\nu)
\xi
(T)}{\nu\alpha}, \frac{\nu\alpha}{2 \beta(2-\nu) \mu(\Omega)}
]$. In particular the assumptions on the sample complexity, the minimum
nonzero singular value of $\La$, and the minimum magnitude nonzero
entry of $\Sa$ are governed by the lower end of this range for $\gamma
$. These assumptions can be weakened if we only require consistency for
a smaller range of values of $\gamma$. The next result conveys this
point with a specific example.\vspace*{-3pt}
%
\begin{corl} \label{theomaincorl}
Consider the same setup and notation as in Theorem~\ref{theomain}.
Suppose that the quantities $\mu(\Omega)$ and $\xi(T)$ satisfy the
assumption of Proposition~\ref{theoirr} for identifiability, and that
$\gamma= \frac{\nu\alpha}{2 \beta(2-\nu) \mu(\Omega)}$ (the upper
end of the range specified in Proposition~\ref{theoirr}), that is,
$\gamma\asymp\frac{1}{\mu(\Omega)}$. Further suppose that: $(1)$~$n
\gtrsim\mu(\Omega)^4  p$; $(2)$ $\lambda_n \asymp\mu(\Omega)
\sqrt
{\frac{p}{n}}$; $(3)$ $\sigma\gtrsim\frac{\mu(\Omega)^2}{\xi(T)}
\sqrt
{\frac{p}{n}}$; $(4)$ $\theta\gtrsim\sqrt{\frac{p}{n}}$. Then with
probability greater than $1 - 2\exp\{-p\}$ we have estimates $(\hat
{S}_n,\hat{L}_n)$ that are algebraically correct, and with the error
bounded as $\g(\hat{S}_n-\Sa, \hat{L}_n - \La) \lesssim\mu
(\Omega)
\sqrt{\frac{p}{n}}$.\vspace*{-3pt}
%
%
%
%
\end{corl}

The proof of this corollary\footnote{By making stronger\vspace*{-1pt}
assumptions on the Fisher information matrix $\F^\ast$, one can further
remove the factor of $\xi(T)$ in the lower bound for $\sigma$.
Specifically, the lower bound $\sigma\gtrsim\mu(\Omega)^3
\sqrt{\frac{p}{n}}$ suffices for consistent estimation if the bounds
defined by the quantities $\alpha_T,\beta_T,\delta_T$ can be
strengthened as described in the remarks at the end of Section
\ref{subsecfi}.} is analogous to that of Theorem~\ref{theomain}. We
emphasize that in practice it is often beneficial to have consistent
estimates for a range of values of $\gamma $ (as in Theorem
\ref{theomain}). Specifically, the stability of the sparsity pattern
and rank of the estimates $(\hat{S}_n,\hat{L}_n)$ for a range of
trade-off parameters is useful in order to choose a suitable value of
$\gamma$, as prior information about the quantities $\mu (\Omega
(\Sa))$ and $\xi(T(\La))$ is not typically available (see Section \ref
{secsims}).\vadjust{\goodbreak}

We remark here that the identifiability conditions of Proposition \ref
{theoirr} are the main sufficient conditions required for Theorem \ref
{theomain} and Corollary~\ref{theomaincorl} to hold. It would be
interesting to obtain necessary conditions as well for these results,
analogous to the necessity and sufficiency of the irrepresentability
conditions for the Lasso \mbox{\cite{ZhaY2006,Wai2009}}.


\subsection{Scaling regimes}
\label{subsecscal}

Next we consider classes of latent-variable models that satisfy the
conditions of Theorem~\ref{theomain}. Recall from Section \ref
{subsecov} that $\mu(\Omega(\Sa)) \leq\degr(\Sa)$. Throughout this
section, we consider latent-variable models in which the low-rank
matrix $\La$ is almost maximally incoherent, that is, $\xi(T(\La))
\sim
\sqrt{\frac{h}{p}}$ so the effect of marginalization over the latent
variables is diffuse across almost all the observed variables. We
suppress the dependence on the quantities $\alpha,\beta,\nu,\psi$
defined in Section~\ref{subsecfi} in our scaling results, and
specifically focus on the trade-off between $\xi(T(\La))$ and $\mu
(\Omega(\Sa))$ for consistent estimation (we also suppress the
dependence of these quantities on $n$). Thus, based on Proposition \ref
{theoirr} we study latent-variable models in which
\[
\xi(T(\La))  \mu(\Omega(\Sa)) = \mathcal{O}\Biggl(\sqrt{\frac
{h}{p}}
\degr(\Sa)\Biggr) = \mathcal{O}(1).
\]
As we describe next, there are nontrivial classes of latent-variable
graphical models in which this condition holds.

\textit{Bounded degree}: The first class of latent-variable models that
we consider are those in which the conditional graphical model among
the observed variables (given by $K^\ast_O$) has constant degree:
\[
\degr(\Sa) = \mathcal{O}(1), \qquad h \sim p.
\]
Such models can be estimated consistently from $n \sim p$ samples. Thus
consistent latent-variable model selection is possible even when the
number of samples and the number of latent variables are on the same
order as the number of observed variables.

\textit{Polylogarithmic degree}: The next class of models that we
consider are those in which the degree of the conditional graphical
model of the observed variables grows polylogarithmically with $p$:
\[
\degr(\Sa) \sim\log(p)^q, \qquad h \sim\frac{p}{\log(p)^{2q}}.
\]
Such latent-variable graphical models can be consistently estimated as
long as $n \sim p $ poly$\log(p)$.

For standard graphical model selection with no latent variables, $\ell
_1$-regular\-ized maximum-likelihood is shown to be consistent with $n =
\mathcal{O}(\log p)$ samples~\cite{RavWRY2008}. On the other hand, our
results prove consistency in the setting with latent variables when $n
= \mathcal{O}(p)$ samples. It would be interesting to study whether
these rates are inherent to latent-variable model selection.

\subsection{Rates for covariance matrix estimation}
\label{subseccov}

Theorem~\ref{theomain} gives conditions under which we can
consistently estimate the sparse and low-rank parts that compose the
marginal concentration matrix $\tilde{K}^\ast_O$. Here we state a
corollary that gives rates for covariance matrix estimation, that is,
the quality of the estimate $(\hat{S}_n-\hat{L}_n)^{-1}$ with respect
to the ``true'' marginal covariance matrix $\Sigma^\ast_O$.\vspace*{-3pt}
%
\begin{corl} \label{theocovcorl}
Under the same conditions as in Theorem~\ref{theomain}, we have with
probability greater than $1 - 2 \exp\{-p\}$ that
\[
\g\bigl(\add^\dag[(\hat{S}_n-\hat{L}_n)^{-1} - \Sigma^\ast_O]\bigr) \leq
\lambda
_n \biggl[1 + \frac{\nu}{6(2-\nu)}\biggr].\vspace*{-2pt}
\]
\end{corl}

This corollary implies that $\|(\hat{S}_n-\hat{L}_n)^{-1} - \Sigma
^\ast
_O\|_2 \lesssim\frac{1}{\xi(T)}\sqrt{\frac{p}{n}}$ based on the
choice of $\lambda_n$ in Theorem~\ref{theomain}, and that $\|(\hat
{S}_n-\hat{L}_n)^{-1} - \Sigma^\ast_O\|_2 \lesssim\mu(\Omega
)\sqrt
{\frac{p}{n}}$ based on the choice of $\lambda_n$ in Corollary \ref
{theomaincorl}.\vspace*{-2pt}


\section{Proofs}\vspace*{-2pt}
\label{secproofs}

\subsection{\texorpdfstring{Proofs of Section \protect\ref{seciden}}{Proofs of Section 3}}
\label{subsecsec3proofs}

Here we give proofs of the results stated in Section~\ref{seciden}.\vspace*{-3pt}
\begin{pf*}{Proof of Lemma~\ref{theorhotspace}}
Since $\rho(T_1,T_2) < 1$, the largest principal angle between $T_1$
and $T_2$ is strictly less than $\frac{\pi}{2}$. Consequently, the
mapping $\proj_{T_2}\dvtx T_1 \rightarrow T_2$ restricted to $T_1$ is
bijective (as it is injective, and the spaces $T_1,T_2$ have the same
dimension). Consider the maximum and minimum gains of $\proj_{T_2}$
restricted to $T_1$; for any $M \in T_1, \|M\|_2 = 1$:
\[
\|\proj_{T_2}(M)\|_2 = \|M + [\proj_{T_2}-\proj_{T_1}](M)\|_2  \in
[1-\rho(T_1,T_2), 1+\rho(T_1,T_2)].
\]
Therefore, we can rewrite $\xi(T_2)$ as follows:
\begin{eqnarray*}
\xi(T_2) &=& {\max_{N \in T_2, \|N\|_2 \leq1}}  \|N\|_\infty = {\max
_{N \in T_2, \|N\|_2 \leq1} } \|\proj_{T_2}(N)\|_\infty\\ &\leq&
\max
_{N \in T_1, \|N\|_2 \leq{1}/({1-\rho(T_1,T_2)})}  \|\proj
_{T_2}(N)\|_\infty\\ &\leq& \max_{N \in T_1, \|N\|_2 \leq
{1}/({1-\rho(T_1,T_2)})}  \bigl[\|N\|_\infty+ \|[\proj_{T_1} - \proj
_{T_2}](N)\|_\infty\bigr] \\ &\leq& \frac{1}{1  -  \rho(T_1,T_2)}
\Bigl[\xi(T_1)  +  {\max_{N \in T_1, \|N\|_2 \leq1}}  \|[\proj_{T_1} -
\proj_{T_2}](N)\|_\infty\Bigr] \\ &\leq& \frac{1}{1  -  \rho
(T_1,T_2)}  \Bigl[\xi(T_1)  +  {\max_{\|N\|_2 \leq1}}  \|[\proj
_{T_1} - \proj_{T_2}](N)\|_2 \Bigr] \\ &\leq& \frac{1}{1  -
\rho
(T_1,T_2)}  [\xi(T_1)  +  \rho(T_1,T_2) ].
\end{eqnarray*}
This concludes the proof of the lemma.\vadjust{\goodbreak}
\end{pf*}
\begin{pf*}{Proof of Lemma~\ref{theorg1}}
We have that $\add^{\dag} \add (S,L) = (S+L, S+L)$; therefore,
$\proj_{\ts} \add^\dag\add\proj _\ts (S,L) = (S + \proj_{\Omega}(L),
\proj_{T}(S) + L)$. We need to bound $\|
S+\proj_{\Omega}(L)\|_{\infty}$ and $\|\proj_{T}(S)+L\|_2$. First, we
have
\begin{eqnarray*}
\|S+\proj_{\Omega}(L)\|_{\infty} &\in& [\|S\|_{\infty} - \|\proj
_{\Omega
}(L)\|_{\infty}, \|S\|_{\infty} + \|\proj_{\Omega}(L)\|_{\infty}]
\\
&\subseteq& [\|S\|_{\infty} - \|L\|_{\infty}, \|S\|_{\infty} + \|L\|
_{\infty}] \\ &\subseteq& [\gamma- \xi(T), \gamma+ \xi(T)].
\end{eqnarray*}
Similarly, one can check that
\begin{eqnarray*}
\|\proj_{T}(S)+L\|_2 &\in& [-\|\proj_{T}(S)\|_2 + \|L\|_2, \|\proj
_{T}(S)\|_2 + \|L\|_2] \\ &\subseteq& [1 - 2\|S\|_2, 1 + 2 \|S\|_2] \\
&\subseteq& [1 - 2 \gamma\mu(\Omega), 1 + 2 \gamma\mu(\Omega)].
\end{eqnarray*}
These two bounds give us the desired result.
\end{pf*}
\begin{pf*}{Proof of Proposition~\ref{theoirr}}
Before proving the two parts of this proposition we make a simple
observation about $\xi(T')$ using the condition that $\rho(T,T')
\leq\frac{\xi(T)}{2}$ by applying Lemma~\ref{theorhotspace}:
\[
\xi(T') \leq\frac{\xi(T) + \rho(T,T')}{1-\rho(T,T')} \leq\frac
{{3\xi(T)}/{2}}{1 - {\xi(T)}/{2}} \leq3 \xi(T).
\]
Here we used the property that $\xi(T) \leq1$ in obtaining the final
inequality. Consequently, noting that $\gamma\in[\frac{3
\beta
(2-\nu) \xi(T)}{\nu\alpha}, \frac{\nu\alpha}{2 \beta(2-\nu)
\mu
(\Omega)} ]$ implies that
%
\begin{equation}\label{eqirreq1}
\chi(\Omega,T',\gamma) = \max\biggl\{\frac{\xi(T')}{\gamma}, 2
\mu
(\Omega) \gamma\biggr\} \leq\frac{\nu\alpha}{\beta(2-\nu)}.
\end{equation}

\textit{Part $1$}: The proof of this step proceeds in a similar manner
to that of Lemma~\ref{theorg1}. First we have for $S \in\Omega, L
\in T'$ with $\|S\|_\infty= \gamma, \|L\|_2 = 1$:
\[
\|\proj_\Omega\F^\ast(S+L)\|_\infty\geq\|\proj_\Omega\F^\ast
S\|
_\infty- \|\proj_\Omega\F^\ast L\|_\infty\geq\alpha\gamma- \|\F
^\ast L\|_\infty\geq\alpha\gamma- \beta\xi(T').
\]
Next, under the same conditions on $S,L$,
\[
\|\proj_{T'} \F^\ast(S+L)\|_2 \geq\|\proj_{T'} \F^\ast L\|_2 -
\|\proj_{T'} \F^\ast S\|_2 \geq\alpha- 2 \|\F^\ast S\|_2 \geq
\alpha
- 2 \beta\mu(\Omega) \gamma.
\]
Combining these last two bounds with~(\ref{eqirreq1}), we conclude that
\begin{eqnarray*}
&& \min_{(S,L) \in\ts,  \|S\|_\infty= \gamma,  \|L\|_2 = 1}  \g
(\proj_\ts\add^\dag\F^\ast\add\proj_\ts(S,L)) \\
&&\qquad
\geq
\alpha- \beta\max\biggl\{\frac{\xi(T')}{\gamma}, 2 \mu(\Omega)
\gamma
\biggr\} \geq\alpha- \frac{\nu\alpha}{2 - \nu} = \frac{2 \alpha
(1-\nu)}{2-\nu} \geq\frac{\alpha}{2},
\end{eqnarray*}
where the final inequality follows from the assumption that $\nu\in
(0,\frac{1}{2}]$.

\textit{Part $2$}: Note that for $S \in\Omega, L \in T'$ with $\|S\|
_\infty\leq\gamma, \|L\|_2 \leq1$,
\[
\|\proj_{\Omega^\bot} \F^\ast(S+L)\|_\infty\leq\|\proj_{\Omega
^\bot}
\F^\ast S\|_\infty+ \|\proj_{\Omega^\bot} \F^\ast L\|_\infty\leq
\delta\gamma+ \beta\xi(T').
\]
Similarly,
\[
\|\proj_{T'^\bot} \F^\ast(S+L)\|_2 \leq\|\proj_{T'^\bot} \F
^\ast S\|
_2 + \|\proj_{T'^\bot} \F^\ast L\|_2 \leq\beta\gamma\mu(\Omega) +
\delta.
\]
Combining these last two bounds with the bounds from the first part, we
have that
\begin{eqnarray*}
&& \|\proj_{\ts^\bot} \add^\dag\F^\ast\add\proj_\ts
(\proj
_\ts\add^\dag\F^\ast\add\proj_\ts)^{-1} \|_{\g
\rightarrow\g}\\
&&\qquad\leq \frac{\delta+ \beta
\max\{{\xi(T')}/{\gamma}, 2 \mu(\Omega) \gamma\}
}{\alpha- \beta\max\{{\xi(T')}/{\gamma}, 2 \mu(\Omega)
\gamma\}} \leq\frac{\delta+ {\nu\alpha}/({2-\nu
})}{\alpha-
{\nu\alpha}/({2-\nu})}\\
&&\qquad\leq\frac{(1-2\nu) \alpha+ {\nu
\alpha
}/({2-\nu})}{\alpha- {\nu\alpha}/({2-\nu})} = 1-\nu.
\end{eqnarray*}
This concludes the proof of the proposition.
\end{pf*}

\subsection{\texorpdfstring{Proof strategy for Theorem \protect\ref{theomain}}{Proof strategy for Theorem 4.1}}
\label{subsecproof}

Standard results\vspace*{1pt} from convex analysis~\cite{Roc1996} state that $(\hat
{S}_n,\hat{L}_n)$ is a minimum of the convex program~(\ref{eqsdp}) if
the zero matrix belongs to the subdifferential of the objective
function evaluated at $(\hat{S}_n,\hat{L}_n)$ [in addition to $(\hat
{S}_n,\hat{L}_n)$ satisfying the constraints]. Elements of the
subdifferentials with respect to the $\ell_1$ norm and the nuclear norm
at a matrix $M$ have the key property that they decompose with respect
to the tangent spaces $\Omega(M)$ and $T(M)$~\cite{Wat92}. This
decomposition property plays a critical role in our analysis. In
particular it states that the optimality conditions consist of two
parts, one part corresponding to the tangent spaces $\Omega$ and $T$
and another corresponding to the normal spaces $\Omega^\bot$ and
$T^\bot$.

Our analysis proceeds by constructing a primal-dual pair of variables
that certify optimality with respect to~(\ref{eqsdp}). Consider the
optimization problem~(\ref{eqsdp}) with the additional (nonconvex)
constraints that the variable $S$ belongs to the algebraic variety of
sparse matrices and that the variable $L$ belongs to the algebraic
variety of low-rank matrices. While this new optimization problem is
nonconvex, it has a very interesting property. At a\vspace*{1pt} globally optimal
solution (and indeed at any locally optimal solution) $(\tilde
{S},\tilde
{L})$ such that $\tilde{S}$ and $\tilde{L}$ are smooth points of the
algebraic varieties of sparse and low-rank matrices, the first-order
optimality conditions state that the Lagrange multipliers corresponding
to the additional variety constraints must lie in the \textit{normal
spaces} $\Omega(\tilde{S})^\bot$ and $T(\tilde{L})^\bot$. This basic
observation, combined with the decomposition property of the
subdifferentials of the $\ell_1$ and nuclear norms, suggests the
following high-level proof strategy: considering the solution $(\tilde
{S},\tilde{L})$ of the variety-constrained problem, we show under
suitable conditions that the second part of the subgradient optimality
conditions of~(\ref{eqsdp}) (without any\vadjust{\goodbreak} variety constraints)
corresponding to components in the normal spaces $\Omega(\tilde
{S})^\bot
$ and $T(\tilde{L})^\bot$ is also satisfied by $(\tilde{S},\tilde{L})$.
Thus, we show that $(\tilde{S},\tilde{L})$ satisfies the optimality
conditions of the \textit{original convex program}~(\ref{eqsdp}).
Consequently $(\tilde{S},\tilde{L})$ is also the optimum of the convex
program~(\ref{eqsdp}). As this estimate is obtained as the solution to
the problem with the variety constraints, the algebraic correctness of
$(\tilde{S},\tilde{L})$ can be directly concluded. We emphasize here
that the variety-constrained optimization problem is used solely as an
analysis tool in order to prove consistency of the estimates provided
by the convex program~(\ref{eqsdp}). The key technical complication is
that the tangent spaces at $\tilde{L}$ and $\La$ are in general
different. We bound the twisting between these tangent spaces by using
the fact that the minimum nonzero singular value of $\La$ is bounded
away from zero (as assumed in Theorem~\ref{theomain}; see also the
supplement~\cite{supp}).


\subsection{Results proved in supplement}



In this section we give the statements of some results that are proved
in a separate supplement~\cite{supp}. These results are critical to the
proof of our main theorem, but they deal mainly with nonstatistical
aspects such as the curvature of the algebraic variety of low-rank
matrices. Recall that $\Omega= \Omega(\Sa)$ and $T = T(\La)$. We also
refer frequently to the constants defined in Theorem~\ref{theomain}.


As the gradient of the log-determinant function is given by a matrix
inverse, a~key step in analyzing the properties of the convex program
(\ref{eqsdp}) is to show that the change in the inverse of a matrix
due to small perturbations is well-approximated by the first-order term
in the Taylor series expansion. Consider the Taylor series of the
inverse of a matrix:
\[
(M + \Delta)^{-1} = M^{-1} - M^{-1} \Delta M^{-1} + R_{M^{-1}}(\Delta),
\]
where
\[
R_{M^{-1}}(\Delta) = M^{-1} \Biggl[\sum_{k=2}^\infty(-\Delta M^{-1})^k
\Biggr].
\]
This infinite sum converges for $\Delta$ sufficiently small. The
following proposition provides a bound on the second-order term
specialized to our setting:

\begin{prop} \label{theorem}
Suppose that $\gamma$ is in the range given by Proposition~\ref
{theoirr}. Further suppose $\Delta_S \in\Omega$, and let $\g
(\Delta
_S,\Delta_L) \leq\frac{1}{2 C_1}$. Then we have that
\[
\g(\add^\dag R_{\Sigma_O^\ast} (\add(\Delta_S,\Delta_L))) \leq
\frac{2
D \psi C_1^2 \g(\Delta_S,\Delta_L)^2}{\xi(T)}.
\]
\end{prop}

Next we analyze the following convex program subject to certain
additional constraints:
%
\begin{eqnarray}\label{eqsdpts}
&&(\hat{S}_\Omega,\hat{L}_{\tilde{T}}) = \mathop{\arg\min}_{S,L} \tr
[(S-L) \Sigma^n_O] - \log\det(S-L)  +  \lambda_n [\gamma\|S\|
_{1} +
\|L\|_\ast] \nonumber\\[-8pt]\\[-8pt]
&&\eqntext{\mbox{s.t. }  S-L \succ0,  S \in\Omega,  L
\in
\tilde{T},}\vspace*{-2pt}
\end{eqnarray}
for some subspace $\tilde{T}$. Comparing~(\ref{eqsdpts}) with the
convex program~(\ref{eqsdp}), we also do not constrain the variable
$L$ to be positive semidefinite in~(\ref{eqsdpts}) for ease of proof
of the next result (see the supplement~\cite{supp} for more details;
recall that the nuclear norm of a positive-semidefinite matrix is equal
to its trace). We show that if $\tilde{T}$ is any tangent space to the
low-rank matrix variety such that $\rho(T,\tilde{T}) \leq\frac{\xi
(T)}{2}$, then we can bound the error $(\Delta_S,\Delta_L) = (\hat
{S}_\Omega- \Sa,\La- \hat{L}_{\tilde{T}})$. Let $\C_{\tilde{T}} =
\proj_{\tilde{T}^\bot}(\La)$ denote the normal component of the true
low-rank matrix at $\tilde{T}$, and let $E_n = \Sigma^n_O - \Sigma
^\ast
_O$ denote the difference between the true marginal covariance and the
sample covariance. The proof of the following result uses Brouwer's
fixed-point theorem~\cite{OrtR1970}, and is inspired by the proof of a
similar result in~\cite{RavWRY2008} for standard sparse graphical model
recovery without latent variables.
%
\begin{prop} \label{theobfpt}
Let the error $(\Delta_S,\Delta_L)$ in the solution of the convex
program~(\ref{eqsdpts}) [with $\tilde{T}$ such that $\rho(\tilde{T},T)
\leq\frac{\xi(T)}{2}$] be as defined above, and define
\[
r = \max\biggl\{\frac{8}{\alpha} [\g(\add^\dag E_n) + \g
(\add^\dag
\F^\ast\C_{\tilde{T}}) + \lambda_n], \|\C_{T'}\|_2
\biggr\}.
\]
If $r \leq\min\{\frac{1}{4 C_1}, \frac{\alpha\xi(T)}{64 D
\psi
C_1^2}\}$ for $\gamma$ as in Proposition~\ref{theoirr}, then
$\g
(\Delta_S,\Delta_L) \leq2r$.
\end{prop}

Finally we give a proposition that summarizes the algebraic component
of our proof.
%
%
\begin{prop}\label{theomainprop}
Assume that $\gamma$ is in the range specified by Proposition~\ref
{theoirr}, $\sigma\geq\frac{C_L \lambda_n}{\xi(T)^2}$, $\theta
\geq
\frac{C_S \lambda_n}{\mu(\Omega)}$, $\g(\add^\dag E_n) \leq
\frac
{\lambda_n \nu}{6(2-\nu)}$, and that $\lambda_n \leq\frac{3
\alpha
(2-\nu)}{16 (3-\nu)} \min\{\frac{1}{4 C_1}$, $\frac{\alpha
\xi
(T)}{64 D \psi C_1^2}\}$. Then there exists a $T'$ and a
corresponding unique solution $(\hat{S}_\Omega,\hat{L}_{T'})$ of
(\ref{eqsdpts}) with $\tilde{T} = T'$ with the following properties:
\begin{longlist}[(2)]
\item[(1)] $\operatorname{sign}(\hat{S}_\Omega) = \operatorname{sign}(\Sa)$ and
$\operatorname{rank}(\hat{L}_{T'}) = \operatorname{rank}(\La)$, with $\hat{L}_{T'}
\succeq
0$. Further $T(\hat{L}_{T'}) = T'$ and $\rho(T,T') \leq\frac{\xi(T)}{4}$.

\item[(2)] Letting $\mathcal{C}_{T'} = \proj_{T'^\bot}(\La)$ we have
that $\g
(\add^\dag\F^\ast\C_{T'}) \leq\frac{\lambda_n \nu}{6 (2-\nu
)}$, and
that $\|\C_{T'}\|_2 \leq\frac{16 (3-\nu) \lambda_n}{3 \alpha
(2-\nu)}$.
\end{longlist}
Further, if $\g(\add^\dag R_{\Sigma^\ast_O} (\add(\hat{S}_\Omega
- \Sa
,\La- \hat{L}_{T'}))) \leq\frac{\lambda_n \nu}{6 (2-\nu)}$,
then the
tangent space constraints $S \in\Omega, L \in T'$ are \textit{inactive}
in~(\ref{eqsdpts}). Consequently the unique solution of~(\ref{eqsdp})
is $(\hat{S}_n,\hat{L}_n) = (\hat{S}_\Omega,\hat{L}_{T'})$.
\end{prop}

\subsection{Probabilistic analysis}

The results given thus far in this section have been completely
deterministic in nature. Here we present the probabilistic component of
our proof by studying the rate at which the sample covariance matrix\vadjust{\goodbreak}
$\Sigma^n_O$ converges to the true covariance matrix $\Sigma_O^\ast$ in
spectral norm. This result is well known and follows directly from
Theorem II.13 in~\cite{DavS2001}; we mainly discuss it here for
completeness and also to show explicitly the dependence on $\psi= \|
\Sigma_O^\ast\|_2$ defined in~(\ref{eqpsi}). See the supplement
\cite
{supp} for a proof.
%
\begin{lemm} \label{theoproblemma}
Let $\psi= \|\Sigma_O^\ast\|_2$. Given any $\delta> 0$ with $\delta
\leq8 \psi$, let the number of samples $n$ be such that $n \geq
\frac{64 p \psi^2}{\delta^2}$. Then we have that
\[
\Pr[\|\Sigma_O^n  -  \Sigma^\ast_O\|_2 \geq\delta]
\leq
2 \exp\biggl\{-\frac{n \delta^2}{128 \psi^2}\biggr\}.
\]
\end{lemm}



The following corollary relates the number of samples required for an
error bound to hold with probability $1-2 \exp\{-p\}$.
%
\begin{corl} \label{theosamples}
Let $\Sigma_O^n$ be the sample covariance formed from $n$ samples of
the observed variables. Set $\delta_n = \sqrt{\frac{128 p \psi
^2}{n}}$. If $n \geq2p$, then
\[
\Pr[\|\Sigma_O^n  -  \Sigma^\ast_O\|_2 \leq\delta_n
]
\geq1 - 2 \exp\{-p\}.
\]
\end{corl}
\begin{pf}
$\!\!\!$Note that $n \geq2p$ implies that $\delta_n \leq8 \psi $, and apply
Lemma~\ref{theoproblemma}.~%
\end{pf}



\subsection{\texorpdfstring{Proof of Theorem \protect\ref{theomain} and Corollary \protect\ref{theocovcorl}}
{Proof of Theorem 4.1 and Corollary 4.3}}
\label{appfinal}

We first combine the results obtained thus far to prove Theorem \ref
{theomain}. Set $E_n = \Sigma_O^n - \Sigma^\ast_O$, set $\delta_n =
\sqrt{\frac{128 p \psi^2}{n}}$, and then set $\lambda_n = \frac
{6 D
\delta_n (2-\nu)}{\xi(T) \nu}$. This\vspace*{1pt} setting of $\lambda_n$ is
equivalent to the specification in the statement of Theorem~\ref{theomain}.
%
%
%
%
%
\begin{pf*}{Proof of Theorem~\ref{theomain}}
We mainly need to show that the various sufficient conditions of
Proposition~\ref{theomainprop} are satisfied. We condition on the event
that \mbox{$\|E_n\|_2 \leq\delta _n$}, which holds with probability greater
than $1 - 2 \exp\{-p\}$ from Corollary~\ref{theosamples} as $n \geq2p$
by assumption. Based on the bound on $n$, we also have that
\[
\delta_n \leq\xi(T)^2 \biggl[\frac{\alpha\nu}{32 (3-\nu) D}
\min
\biggl\{ \frac{1}{4 C_1},\frac{\alpha\nu}{256 D (3-\nu) \psi C_1^2}
\biggr\} \biggr].
\]
In particular, these bounds imply that
%
\begin{eqnarray} \label{eqdelta}
\delta_n &\leq&\frac{\alpha\xi(T) \nu}{32 (3-\nu) D} \min
\biggl\{ \frac
{1}{4 C_1},\frac{\alpha\xi(T)}{64 D \psi C_1^2}
\biggr\};\nonumber\\[-8pt]\\[-8pt]
\delta_n &\leq&\frac{\alpha^2 \xi(T)^2 \nu^2}{8192 \psi C_1^2
(3-\nu)^2
D^2}.\nonumber
\end{eqnarray}
Both these weaker bounds are used later.


Based on the assumptions of Theorem~\ref{theomain}, the requirements
of Proposition~\ref{theomainprop} on $\sigma$ and $\theta$ are
satisfied. Next we verify the bounds on $\lambda_n$ and $\g(\add
^\dag
E_n)$.\vadjust{\goodbreak} Based on the setting of $\lambda_n$ above and the bound on
$\delta_n$ from~(\ref{eqdelta}), we have that
\[
\lambda_n = \frac{6 D (2-\nu) \delta_n}{\xi(T) \nu} \leq\frac{3
\alpha(2-\nu)}{16(3-\nu)} \min\biggl\{\frac{1}{4C_1},\frac
{\alpha\xi
(T)}{64 D \psi C_1^2} \biggr\}.
\]
Next we combine the facts that $\lambda_n = \frac{6 D \delta_n
(2-\nu
)}{\xi(T) \nu}$ and that $\|E_n\|_2 \leq\delta_n$ to conclude that
%
\begin{equation}\label{eqsamplambound}
\g(\add^\dag E_n) \leq\frac{D \delta_n}{\xi(T)} = \frac
{\lambda_n \nu
}{6 (2 - \nu)}.
\end{equation}

Thus, we have from Proposition~\ref{theomainprop} that there exists a
$T'$ and corresponding solution $(\hat{S}_\Omega,\hat{L}_{T'})$ of
(\ref{eqsdpts}) with the prescribed properties. Next we apply
Proposition~\ref{theobfpt} with $\tilde{T} = T'$ to bound the error
$(\hat{S}_\Omega- \Sa, \La- \hat{L}_{T'})$. Noting that $\rho(T,T')
\leq\frac{\xi(T)}{4}$, we have that
%
\begin{eqnarray}
\frac{8}{\alpha}[\g(\add^\dag E_n) + \g(\add^\dag\F^\ast
\C
_{T'}) + \lambda_n ] &\leq& \frac{8}{\alpha} \biggl[\frac
{\nu}{3
(2-\nu)} + 1 \biggr] \lambda_n \nonumber\\
\label{eqlambound} &=& \frac{16(3-\nu
)\lambda
_n}{3 \alpha(2-\nu)} \\
\label{eqdelbound} &=& \frac{32(3-\nu
)D}{\alpha\xi(T) \nu} \delta_n \\
\label{eqfinalmainineq}
&\leq& \min
\biggl\{\frac{1}{4C_1},\frac{\alpha\xi(T)}{64 D \psi C_1^2} \biggr\}.
\end{eqnarray}
In the first inequality we used the fact that $\g(\add^\dag E_n) \leq
\frac{\lambda_n \nu}{6(2-\nu)}$ (from above) and that $\g(\add
^\dag\F
^\ast\C_{T'})$ is similarly bounded (from Proposition \ref
{theomainprop}). In the second equality we used the relation $\lambda
_n = \frac{6 D \delta_n (2-\nu)}{\xi(T) \nu}$. In the final inequality
we used the bound on $\delta_n$ from~(\ref{eqdelta}). This satisfies
one of the requirements of Proposition~\ref{theobfpt}. The second
requirement of Proposition~\ref{theobfpt} on $\|\C_{T'}\|_2$ is also
similarly satisfied as we have that $\|\C_{T'}\|_2 \leq\frac
{16(3-\nu
)\lambda_n}{3 \alpha(2-\nu)}$ from Proposition~\ref{theomainprop},
and we use the same sequence of inequalities as above. Thus we conclude
from Proposition~\ref{theobfpt} and from~(\ref{eqlambound}) that
%
\begin{equation}\label{eqfinerrbound}
\g(\hat{S}_\Omega- \Sa,\La- \hat{L}_{T'}) \leq\frac{32(3-\nu
)\lambda
_n}{3 \alpha(2-\nu)} \lesssim\frac{1}{\xi(T)} \sqrt{\frac{p}{n}}.
\end{equation}
Here the last inequality follows from the bound on $\lambda_n$.

If we show that $(\hat{S}_n,\hat{L}_n) = (\hat{S}_\Omega, \hat
{L}_{T'})$, we can conclude the proof of Theorem~\ref{theomain} since
algebraic correctness of $(\hat{S}_\Omega,\hat{L}_{T'})$ holds from
Proposition~\ref{theomainprop} and the estimation error bound follows
from~(\ref{eqfinerrbound}). In order to complete this final step, we
again\vspace*{1pt} revert to Proposition~\ref{theomainprop} and prove the requisite
bound on $\g(\add^\dag R_{\Sigma_O^\ast}(\add(\hat{S}_\Omega-
\Sa,\La
- \hat{L}_{T'})))$.

Since the bound~(\ref{eqfinerrbound}) combined with the inequality
(\ref{eqfinalmainineq}) satisfies the condition of Proposition \ref
{theorem} [i.e., we have that $\g(\hat{S}_\Omega- \Sa,\La- \hat
{L}_{T'}) \leq\frac{1}{2C_1}$]:
\begin{eqnarray*}
\g\bigl(\add^\dag R_{\Sigma_O^\ast}\bigl(\add(\hat{S}_\Omega- \Sa,\La-
\hat
{L}_{T'})\bigr)\bigr) &\leq& \frac{2 D \psi C_1^2}{\xi(T)} \g(\hat
{S}_\Omega-
\Sa,\La- \hat{L}_{T'})^2 \\ &\leq& \frac{2 D \psi C_1^2}{\xi(T)}
\biggl(\frac{64(3-\nu)D}{\alpha\xi(T) \nu}\biggr)^2 \delta_n^2 \\ &=&
\biggl[\frac{8192 \psi C_1^2 (3-\nu)^2 D^2}{\alpha^2 \xi(T)^2 \nu^2}
\delta
_n\biggr] \frac{D \delta_n}{\xi(T)} \\ &\leq& \frac{D \delta
_n}{\xi
(T)} \\ &=& \frac{\lambda_n \nu}{6 (2-\nu)}.
\end{eqnarray*}
In the second inequality we used~(\ref{eqdelbound}) and (\ref
{eqfinerrbound}), in the final inequality we used the bound (\ref
{eqdelta}) on $\delta_n$, and in the final equality we used the
relation $\lambda_n = \frac{6 D \delta_n (2-\nu)}{\xi(T) \nu}$.
\end{pf*}
\begin{pf*}{Proof of Corollary~\ref{theocovcorl}}
Based on the optimality conditions of the modified convex program
(\ref{eqsdpts}), we have that
\[
\g\bigl(\add^\dag[(\hat{S}_n-\hat{L}_n)^{-1} - \Sigma^n_O]\bigr) \leq
\lambda_n.
\]
Combining this with the bound~(\ref{eqsamplambound}) yields the
desired result.
\end{pf*}

\section{Simulation results}
\label{secsims}

In this section we give experimental demonstration of the consistency
of our estimator~(\ref{eqsdp}) on synthetic examples, and its
effectiveness in modeling real-world stock return data. Our choices of
$\lambda_n$ and $\gamma$ are guided by Theorem~\ref{theomain}.
Specifically, we choose $\lambda_n$ to be proportional to $\sqrt
{\frac
{p}{n}}$. For $\gamma$ we observe that the support/sign-pattern and the
rank of the solution $(\hat{S}_n,\hat{L}_n)$ are the same for a \textit
{range} of values of $\gamma$. Therefore one could solve the convex
program~(\ref{eqsdp}) for several values of $\gamma$, and choose a
solution in a suitable range in which the sign-pattern and rank of the
solution are stable (see~\cite{ChaSPW2009} for details). In practical
problems with real-world data these parameters may be chosen via
cross-validation (it would be of interest to consider methods such as
those developed in~\cite{MeiB2010}). For small problem instances we
solve the convex program~(\ref{eqsdp}) using a combination of YALMIP
\cite{Lof2004} and SDPT3~\cite{TohTT}. For larger problem instances we
use the special-purpose solver LogdetPPA~\cite{WanST2009} developed for
log-determinant semidefinite programs.

\subsection{Synthetic data}
In the first set of experiments we consider a setting in which we have
access to samples of the observed variables of a latent-variable
graphical model. We consider several latent-variable\vadjust{\goodbreak} Gaussian graphical
models. The first model consists of $p=36$ observed variables and $h =
2$ latent variables. The conditional graphical model structure of the
observed variables is a cycle with the edge partial correlation
coefficients equal to $0.25$; thus, this conditional model is specified
by a sparse graphical model with degree $2$. The second model is the
same as the first one, but with $h = 3$ latent variables. The third
model consists of $h = 1$ latent variable, and the conditional
graphical model structure of the observed variables is given by a $6
\times6$ nearest-neighbor grid (i.e., $p=36$ and degree $4$) with the
partial correlation coefficients of the edges equal to $0.15$. In all
three of these models each latent variable is connected to a random
subset of $80\%$ of the observed variables (and the partial correlation
coefficients corresponding to these edges are also random). Therefore
the effect of the latent variables is ``spread out'' over most of the
observed variables, that is, the low-rank matrix summarizing the effect
of the latent variables is incoherent.

\begin{figure}

\includegraphics{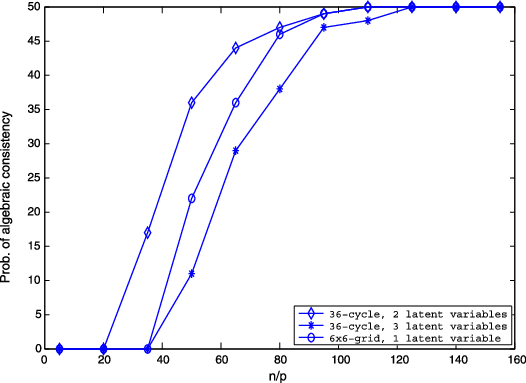}

\caption{Synthetic data: plot showing probability of algebraically
correct estimation. The three models studied are \textup{(a)} 36-node
conditional graphical model given by a cycle with $h=2$ latent
variables, \textup{(b)} 36-node conditional graphical model given by a cycle
with $h=3$ latent variables and \textup{(c)} 36-node conditional graphical
model given by a $6 \times6$ grid with $h = 1$ latent variable. For
each plotted point, the probability of algebraically correct estimation
is obtained over $50$ random trials.} \label{figfig1}
\end{figure}

For each model we generate $n$ samples of the observed variables, and
use the resulting sample covariance $\Sigma_O^n$ as input to our convex
program~(\ref{eqsdp}). Figure~\ref{figfig1} shows the probability of
obtaining algebraically correct estimates as a function of~$n$. This
probability is evaluated over $50$ experiments for each value of $n$.
In all of these cases standard graphical model selection applied
directly to the observed variables is not useful as the marginal
concentration matrix of the observed variables is not well-approximated
by a sparse matrix. These experiments agree with our theoretical
results that the convex program~(\ref{eqsdp}) is an algebraically
consistent estimator of a latent-variable model given (sufficiently
many) samples of only the observed variables.\vspace*{-3pt}


\subsection{Stock return data}
In the next experiment we model the statistical structure of monthly
stock returns of 84 companies in the S\&P 100 index from 1990 to 2007;
we disregard 16 companies that were listed after 1990. The number of
samples $n$ is equal to $216$. We compute the sample covariance based
on these returns and use this as input to~(\ref{eqsdp}).

\begin{figure}

\includegraphics{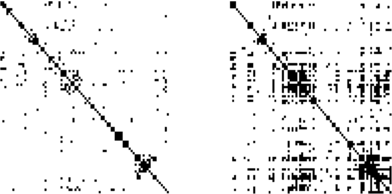}

\caption{Stock returns: the figure on the left shows the sparsity
pattern (black denotes an edge, and white denotes no edge) of the
concentration matrix of the conditional graphical model (135 edges) of
the stock returns, conditioned on five latent variables, in a
latent-variable graphical model (total number of parameters equals
$639$). This model is learned using (\protect\ref{eqsdp}), and the KL
divergence with respect to a Gaussian distribution specified by the
sample covariance is $17.7$. The figure on the right shows the
concentration matrix of the graphical model (646 edges) of the stock
returns, learned using standard sparse graphical model selection based
on solving an $\ell_1$-regularized maximum-likelihood program (total
number of parameters equals $730$). The KL divergence between this
distribution and a Gaussian distribution specified by the sample
covariance is $44.4$.} \label{figfig2}\vspace*{-3pt}
\end{figure}

The model learned using~(\ref{eqsdp}) for suitable values of $\lambda
_n,\gamma$ consists of $h = 5$ latent variables, and the conditional
graphical model structure of the stock returns conditioned on these
latent components consists of $135$ edges. Therefore the number of
parameters in the model is $84 + 135 + (5 \times84) = 639$. The
resulting KL divergence between the distribution specified by this
model and a Gaussian distribution specified by the sample covariance is
$17.7$. Figure~\ref{figfig2} (left) shows the \textit{conditional}
graphical model structure. The strongest edges in this conditional
graphical model, as measured by partial correlation, are between Baker
Hughes--Schlumberger, A.T.\&T.--Verizon, Merrill Lynch--Morgan
Stanley, Halliburton--Baker Hughes, Intel--Texas Instruments,
Apple--Dell, and Microsoft--Dell. It is of interest to note that in the
Standard Industrial Classification\footnote{See the U.S. SEC website at
\url{http://www.sec.gov/info/edgar/siccodes.htm}.} system for grouping these
companies, several of these pairs are in different classes. As
mentioned in Section~\ref{subsecps}, our method estimates a low-rank
matrix that summarizes the effect of the latent variables; in order to
factorize this low-rank matrix, for example, into sparse factors, one
could use methods such as those described in~\cite{WitTH2009}.\vadjust{\goodbreak}

We compare these results to those obtained using a sparse graphical
model learned using $\ell_1$-regularized maximum-likelihood (see, e.g.,
\cite{RavWRY2008}), without introducing any latent variables.
Figure~\ref{figfig2} (right) shows this graphical model structure. The
number of edges in this model is $646$ (the total number of parameters
is equal to $646 + 84 = 730$), and the resulting KL divergence between
this distribution and a Gaussian distribution specified by the sample
covariance is $44.4$.


These results suggest that a latent-variable graphical model is better
suited than a standard sparse graphical model for modeling stock
returns. This is likely due to the presence of global, long-range
correlations in stock return data that are better modeled via latent variables.

\section{Discussion}
\label{secconc}

We have studied the problem of modeling the statistical structure of a
collection of random variables as a sparse graphical model conditioned
on a few additional latent components. As a first contribution we
described conditions under which such latent-variable graphical models
are identifiable given samples of only the observed variables. We also
proposed a convex program based on $\ell_1$ and nuclear norm
regularized maximum-likelihood for latent-variable graphical model
selection. Given samples of the observed variables of a latent-variable
Gaussian model, we proved that this convex program provides consistent
estimates of the number of latent components as well as the conditional
graphical model structure among the observed variables conditioned on
the latent components. Our analysis holds in the high-dimensional
regime in which the number of observed/latent variables are allowed to
grow with the number of samples of the observed variables. These
theoretical predictions are verified via a set of experiments on
synthetic data. We also demonstrate the effectiveness of our approach
in modeling real-world stock return data.

Several questions arise that are worthy of further investigation. While~(\ref{eqsdp})
can be solved in polynomial time using off-the-shelf
solvers, it is preferable to develop more efficient special-purpose
solvers to scale to massive datasets by taking advantage of the
structure of~(\ref{eqsdp}). It is also of interest to develop
statistically consistent convex optimization methods for
latent-variable modeling with non-Gaussian variables, for example, for
categorical data.

\section*{Acknowledgments}

We would like to thank James Saunderson and Myung Jin Choi for helpful
discussions, and Kim-Chuan Toh for kindly providing us specialized code
to solve larger instances of our convex program.

\begin{supplement}[id=suppA]
\stitle{Supplement to ``Latent variable
graphical model selection via convex optimization''}
\slink[doi]{10.1214/11-AOS949SUPP} 
\sdatatype{.pdf}
\sfilename{aos949\_supp.pdf}
\sdescription{Due to space constraints, we have
moved some technical proofs to a supplementary document
\cite{supp}.}
\end{supplement}


\printaddresses

\end{document}